\documentclass[a4paper]{article}
 \usepackage{xcolor}
 \colorlet{RED}{red}
 \usepackage{graphics}
 \usepackage{tikz}
 \usepackage{graphicx}
 \usepackage{subfigure}
 \usepackage{epsfig}
\usepackage{authblk}

 \usepackage{amsthm}
 \usepackage{mathrsfs}
\usepackage{epsfig,graphicx,amssymb,amsmath,cite}
\usepackage{ulem}
\usepackage{amsfonts}
\usepackage{savesym}
\usepackage{listings}
\usepackage{titlesec}
\usepackage{bm}
\usepackage{geometry}
\usepackage{appendix}
\usepackage{epstopdf}
\usepackage{ulem}

\newtheorem{thm}{\textbf{Theorem}}

\newtheorem{lem}{\textbf{Lemma}}

\newtheorem{remark}{\textbf{Remark}}

\newcommand{\argmin}{\mathop{\arg\min}}
\newcommand{\argmax}{\mathop{\arg\max}}

\newtheorem{proposal}{\textbf{Proposal}}

\newcommand{\TT}{{\mathrm{T}}}
\newcommand{\dd}{{\mathrm{d}}}

\usetikzlibrary{arrows,automata}

\providecommand{\keywords}[1]{\textbf{\textit{Keywords:}} #1}

\title{\bf \Large On the Mathematics of RNA Velocity II: Algorithmic Aspects}
\author[1,2,3]{Tiejun Li\thanks{tieli@pku.edu.cn}}
\author[3]{Yizhuo Wang\thanks{jiuqie@pku.edu.cn}}
\author[1]{Guoguo Yang\thanks{ygj512@hotmail.com}}
\author[4]{Peijie Zhou\thanks{peijiez1@uci.edu}}

\affil[1]{\small LMAM and School of Mathematical Sciences, Peking University, Beijing 100871, China}
\affil[2]{Center for Machine Learning Research, Peking University, Beijing 100871, China}
\affil[3]{Center for Data Science, Peking University, Beijing 100871, China}
\affil[4]{Department of Mathematics, University of California, Irvine, Irvine, CA 92697, USA}

\date{}

\begin{document}

\maketitle

\begin{abstract}
In a previous paper [CSIAM Trans. Appl. Math. 2 (2021), 1-55], the authors proposed a  theoretical framework for the analysis of RNA velocity, which is a promising concept in scRNA-seq data analysis to reveal the cell state-transition dynamical processes underlying snapshot data.  The current paper is devoted to the algorithmic study of some key components in RNA velocity workflow.  Four important points are addressed in this paper: (1) We construct a rational time-scale fixation method which can determine the global gene-shared latent time for cells. (2) We present an uncertainty quantification strategy for the inferred parameters obtained through the EM algorithm. (3) We establish the optimal criterion for the choice of velocity kernel bandwidth with respect to the sample size in the downstream analysis and discuss its implications. (4) We propose a temporal distance estimation approach between two cell clusters along the cellular development path.  Some illustrative numerical tests are also carried out to verify our analysis. These results are intended to provide tools and insights in further development of RNA velocity type methods in the future.
\end{abstract}

\section{Introduction}

\begin{figure}[h]
    \centering
    \includegraphics{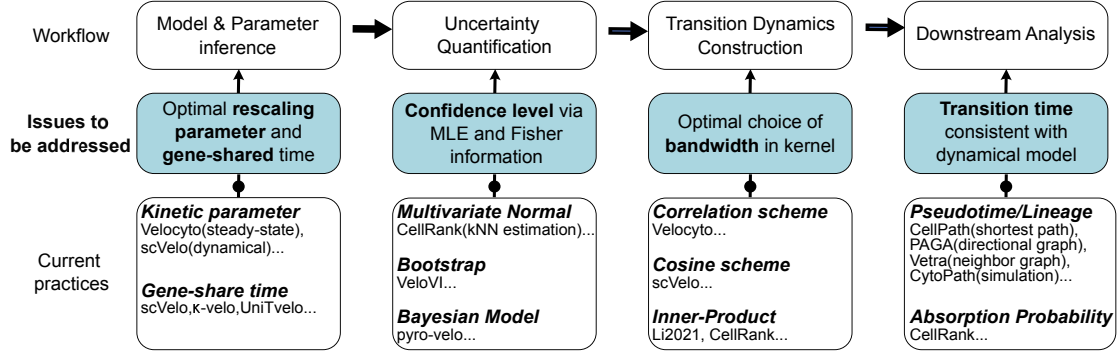}
    \caption{The computational workflow of RNA velocity analysis and under-addressed issues.}
    \label{fig:intro}
\end{figure}

The development of single-cell RNA sequencing (scRNA-seq) technology has revolutionized the  resolution and capability to dissect the cell-fate determination process\cite{Tang2009}. However, traditional scRNA-seq datasets only provide static snapshots of gene expression among cells at a certain time point, which lack the direct temporal information to infer the dynamics of cell state transitions\cite{weinreb2018fundamental}. To address this limitation, the RNA velocity method \cite{MannoRNA2018} utilizes both unspliced and spliced counts in scRNA-seq data to model and infer the dynamics of mRNA expression and splicing process, allowing the prediction of gene expression changes over time, and the specification of directionality during development. The method has been applied widely in different biological systems
\cite{gorin2022rna,Farrell2022,atta2022veloviz}, and the computational workflow of RNA velocity analysis has been established and undergone rapid development\cite{Zhang2021,MannoRNA2018,BergenGene2020,Lange2022} (Fig. \ref{fig:intro}).

To improve the effectiveness and robustness of RNA velocity analysis, various algorithmic modifications have been proposed throughout the computational workflow. For the parameter inference step, scVelo\cite{BergenGene2020} utilizes an Expectation-Maximization (EM) procedure between latent time specification and kinetic parameter update to generalize the steady-state assumption to the transient dynamical process. In addition, $\kappa$-velo\cite{marot2022towards} proposes to calculate a gene-shared latent time for each cell by approximating the traveling time with the number of cells in-between, and UniTvelo\cite{gao2022unitvelo} calculates the unified latent time by aggregating the gene-specific time quantiles. Recently, VeloVAE utilizes variational Bayesian inference and autoencoder to compute the gene-shared latent time and cell latent state \cite{gu2022bayesian}. To account for the uncertainty of inferred parameters incurred by noise and sparsity in spliced or unspliced counts, CellRank \cite{Lange2022} adopts the multivariate normal model to quantify the velocity distribution, while VeloVI \cite{gayoso2022deep} employs the bootstrap strategy. Recently, pyro-Velo \cite{qin2022pyro} proposes a Bayesian approach to model the posterior distribution of parameters.

Based on the inferred RNA velocity, downstream dynamical analysis tools such as low-dimensional embedding\cite{qiao2021representation,atta2022veloviz} and trajectory inference\cite{Zhang2021,gao2022unitvelo,liu2022dynamical} are developed by leveraging the cell-cell neighbor graph directed by the velocities. Pertinent to such methods is the construction of a cellular random walk transition probability (or weight) matrix, which is induced by the velocity-based data kernels. Our previous theoretical work \cite{Li1} has elucidated that  different choices of velocity kernels can result in various forms of differential equations as the continuum limit. The resulting Markov chain of cell-state transitions can provide further insights into the underlying dynamics, including the transition paths to quantify development routes \cite{Li1,zhou2021dissecting}, and the absorption probabilities to quantify cell-fate commitment likelihood \cite{setty2019characterization,Lange2022}.

Despite the success of algorithm developments, deeper analysis is still necessary to understand the rationales of the algorithm design and address the unresolved issues. For instance, in parameter estimation, the optimal choice of re-scaling parameter and therefore the determination of gene-shared latent time beyond heuristic strategies have yet to be determined. In fact, extracting the hidden time information from the snapshot data is still a central issue in scRNA-seq data analysis. Meanwhile, the uncertainty quantification of parameters could benefit from a rigorous confidence level analysis of the EM algorithm. In the construction of random walks, the appropriate choice of kernel bandwidth relies on the numerical analysis of convergence order to continuum limit. In addition, simulations and benchmarks are important to validate the statistical and numerical analysis results.

As a continuation of our previous work \cite{Li1}, in this paper we will study the mathematics of RNA velocity analysis from the  perspective of algorithms. The main contributions of this paper toward the current computation workflow could be summarized as follows:

\begin{itemize}
    \item \textbf{Parameter Inference.} To determine the gene-latent time \cite{BergenGene2020}, we formulate an optimization framework to determine the gene-specific rescaling parameters and propose the numerical scheme to efficiently tackle the considered problem.
    \item \textbf{Uncertainty Quantification.} Employing the asymptotic theory for EM algorithm\cite{MengUsing1991}, we rigorously derive the confidence level for the kinetic parameters in the dynamical RNA velocity model.
    \item \textbf{Random Walk Construction.} Addressing the finite sample-size issue in  computation, we analyze the variance and bias of the approximation to continuum equation, and find an optimal criterion to determine kernel bandwidth and sample size for the velocity kernel which induces the cellular random walk dynamics.
    \item \textbf{Downstream Analysis.} To perform lineage inference that is consistent with RNA velocity dynamics, we propose the first hitting time analysis to quantify the transition time.
\end{itemize}

The rest of this paper is organized as follows. According to the contents stated above, we have studied them in Sections 2, 3, 4 and 5, respectively, and given corresponding numerical illustrations in each section. Finally, we make the conclusion. Some proof details are left in the Appendix.

\section{Fixation of time re-scaling constants}

The current RNA velocity models \cite{MannoRNA2018,BergenGene2020} assume that the genes are independent, thus in dynamical parameter inference\cite{BergenGene2020}, the latent times of the cells are obtained for each gene independently, and we lack a determination step of time-rescaling factors to form a globally consistent gene-shared time. To rationally infer the global cell hidden times as well as determine the rescaling factors, we propose an optimization framework to address this issue. It could be served as a good starting point for more delicate inference on the latent time by taking into account more technical details.

\subsection{Problem setup}

Suppose that in a considered scRNA-seq measurement, we have $d$ genes with the label $g=1,2.\dots,d$ and $n$ cells with the label $c=1,2,\ldots,n$. Similar to the previous work, we utilize the deterministic dynamical model
\begin{equation}\label{deter_model}
\begin{aligned}
\frac{\mathrm{d} u}{\mathrm{~d} t} &=\alpha^{\text {on}/\text {off}}(t)-\beta u(t), \\
\frac{\mathrm{d} s}{\mathrm{~d} t} &=\beta u(t)-\gamma s(t),
\end{aligned}
\end{equation}
to describe the transcriptional process of each gene, and the individual genes are independent of each other.
Here $t \geq 0,\left.(u(t), s(t))\right|_{t=0}=\left(u_{0}, s_{0}\right)$, and
$$
\alpha^{\text {on}/\text {off}}(t)= \begin{cases}\alpha^{\text {on}} ,& t \leq t_{s}, \\ \alpha^{\text {off}}=0, & t>t_{s},\end{cases}
$$
where $t_{s}$ is the switching time of the transcriptional process at which transcription rate $\alpha$ turns to $0$. The variables $u(t)$ and $s(t)$ are the abundance of unspliced and spliced mRNA in the cell measured at time $t$, respectively. In general, the resulting data are not time-resolved and $t $ is a latent variable. Likewise, the transcriptional state of the cell (on/off) is an unknown variable, and the rates $\alpha^{\text {on}}$, $\beta$, and $\gamma$ cannot be directly measured experimentally.

In the inference process, we need to solve the equation and infer the kinetics of splicing controlled by parameters: transcription rate $\alpha^{\text {on}}$, splicing rate $\beta$ and degradation rate $\gamma$; latent variable time $t$.  We usually infer the parameters for each gene separately under the independent gene assumption, which leaves the relative size of the parameters of  genes as an unsolved problem. As the system has the following scale invariance property \cite{Li1}, i.e., if we define the parameter $\theta=(\theta_r,t_s)$, in which $\theta_r=(\alpha,\beta,\gamma)$, then the following equation holds:
\begin{equation}
    \big(u(t;\theta_r,t_s),s(t;\theta_r,t_s)\big)=\big(u(\kappa t;\theta_r/\kappa,\kappa t_s),s(\kappa t;\theta_r/\kappa,\kappa t_s)\big),
\end{equation}
where $\kappa>0$ is the scaling parameter. In the inference we usually keep $\beta_g=1$ at first while optimizing other parameters for each specific gene $g$, which essentially infers $\alpha_g/\beta_g$ and $\gamma_g/\beta_g$ due to the scale invariance. When considering the high-dimensional velocity and the corresponding low-dimensional projection, the scale needs to be adjusted among all genes, that is, the scaling parameters $(\beta_g)_g$ need to be determined for each gene. As this parameter appears in the final RNA velocity, its choice will highly affect the lineage inference in downstream analysis. Also, computing the gene-latent time required us to find out the scaling parameters. Indeed, this is an important under-addressed issue in scRNA-seq data analysis.

Assume that we have already inferred the unscaled parameters $\alpha_g$, $\beta_g=1$, $\gamma_g$  for each gene, with the gene-specific cell time matrix $T=(t_{cg})\in(\mathbb{R}^+\cup\{0\})^{n\times d}$. Our goal is to infer the gene-shared latent time $t_c$ for each cell as well as determine the rescaling parameters $\beta_g$ for different genes. Below we will propose two optimization approaches to tackle this issue.

\subsection{Gene-shared time through optimization}\label{sec:2.2}

To obtain the gene-shared latent time in any given cell, we reason it to be as consistent as possible with the respective rescaled time for each gene within the cell.
Denote by $\beta=(\beta_g)$ or $x=(x_g)=(\beta_g^{-1})\in \mathbb{R}^d$ the time re-scaling parameters for the genes, and $t=(t_c)\in \mathbb{R}^n$ the gene-shared latent time for cells to be optimized. We formulate the above consistency intuition through two proposals.

Our first proposal is based on the model
\begin{equation}\label{eq:prop1}
t_{cg}\beta_g^{-1}=t_c+\epsilon_{cg},\quad \epsilon_{cg}\sim N(0,\sigma^2) \quad \textrm{for } c=1,\ldots,n;\ g=1,\ldots,d.
\end{equation}
Here $t_{cg}$ is the inferred gene-specific time with $\beta_g=1$, and $t_{cg}\beta_g^{-1}$ is the rescaled time with $\beta_g$, and \eqref{eq:prop1} means that the rescaled time should be consistent with a global gene-shared common time $t_c$ upon removing some noise. With this setup, we can determine $x$ and $t$ with the following formulation.

\begin{proposal}[Inference with Multiplicative Noise] The gene-shared latent time $t$ and rescaling parameters $x$ can be determined through the minimization problem
\begin{equation}\label{eq:timescale}
    (x^*, t^*) = \argmin_{\Vert t\Vert=1; x\succ 0, t\succeq 0}\Vert TX-t\bm{1}^{\mathrm{T}}\Vert_F^2,
\end{equation}
where $x\succ 0$, $t\succeq 0$ means that $x,t$ have positive or non-negative components, respectively, $\bm{1}=(1,1,\ldots,1)^{\mathrm{T}}\in\mathbb{R}^d$, $X=\mathrm{diag}(x_1,\ldots,x_d)\in\mathbb{R}^{d\times d}$ is the diagonal matrix formed by the components of $x$, and
$$\Vert A \Vert_F:=\Big(\sum_{ij} a^2_{ij}\Big)^\frac12 =({\rm tr}(AA^{\mathrm{T}}))^\frac12 \quad \textrm{ for } A=(a_{ij})$$
denotes the Frobenius norm ($F$-norm) of a matrix.
\end{proposal}

\begin{thm}\label{thm:1} Assume that the inferred gene-specific cell time matrix $T$ satisfies the condition
\begin{equation}\label{eq:Irreduc}
T\in \mathcal{T}=\big\{T\in(\mathbb{R}^+\cup\{0\})^{n\times d}\ |\ T^{\mathrm{T}}T\text{ is irreducible}\big\}.
\end{equation}
Then the optimization problem \eqref{eq:timescale} has the unique solution $x^*=d \lambda_1^{-1/2} W v_1$, where $v_1$ is the $\ell^2$-unit eigenvector corresponding to the maximal eigenvalue $\lambda_1$ of
\begin{equation}
H=W^\mathrm{T}T^\mathrm{T}TW, \quad \textrm{ where } W:=\mathrm{diag}(w_1,\ldots,w_d),\ \ w_g=1/\Vert t_{\bullet g} \Vert,\  g=1,\ldots,d
\end{equation}
and it has positive components. The global gene-shared common time
$$t^*=TWv_1/\|TWv_1\|.$$
\end{thm}
\begin{proof}
To solve the problem \eqref{eq:timescale}, we note that when $x$ is fixed, the optimization
$$
\min_{t\succeq 0}\Vert TX-t\bm{1}^{\mathrm{T}}\Vert_F^2
$$
turns out to be a least squares problem, and the minimum point is $t = Tx/d.$
Substituting it back to \eqref{eq:timescale}, we get
$$
\begin{aligned}
x^*&=\argmin_{\Vert Tx\Vert=d,x\succ 0}\Big\| TX-\frac1d TXE\Big\|_F^2.
\end{aligned}
$$
where the matrix $E:=\bm{1}\bm{1}^{\mathrm{T}} \in \mathbb{R}^{d\times d}$. Define $C=I_{d}-E$, which satisfies $C^{\mathrm{T}}=C$ and $C^2=C$. We have
$$
\begin{aligned}
x^*&=\argmin_{\Vert Tx\Vert=d,x\succ 0}\Vert TXC\Vert_F^2.
\end{aligned}
$$
By the definition of the $F$-norm, we have
\begin{equation}\label{eq:5}
x^*=\argmin_{\Vert Tx\Vert=d,x\succ 0}\text{tr}(TXCC^{\mathrm{T}}X^{\mathrm{T}}T^{\mathrm{T}}).
\end{equation}
Denote by $A\circ B$ the Hadamard product of matrices $A$ and $B$ defined as $A\circ B=(a_{ij}b_{ij})$ for $A=(a_{ij})$ and $B=(b_{ij})$. It is not difficult to find that \eqref{eq:5} is equivalent to
\begin{equation}\label{eq:4}
x^*=\argmin_{\Vert Tx\Vert=d,x\succ 0}x^{\mathrm{T}}Mx,
\end{equation}
 where $M=(T^{\mathrm{T}} T)\circ C$.

Denote by $t_{\bullet g}$ the vector formed by $(t_{cg})_c$ for a fixed gene $g$. By the irreducibility condition \eqref{eq:Irreduc},
 we have $\|t_{\bullet g}\|>0$ for any $g$, thus $W$ is well-defined. Further
 note that $M=(T^\TT T)\circ C=W^{-2} - T^\TT T/d$, the problem \eqref{eq:4} is equivalent to
\begin{equation}\label{eq:10}
x^*=\argmin_{\Vert Tx\Vert=d,x\succ 0}x^{\mathrm{T}}W^{-2}x.
\end{equation}
Suppose
$$H=Q^\TT \Lambda Q,\quad Q^\TT=(v_1,v_2,\ldots,v_d)$$
where $Q^\TT Q= I_d$, $\Lambda={\rm diag}(\lambda_1,\ldots,\lambda_d)$ and $\lambda_1\ge \lambda_2\ge \cdots\ge \lambda_d\ge 0$. We have $Hv_k=\lambda_k v_k$ for $k=1,\ldots,d$. Ignoring the positivity constraint $x\succ 0$, we can find that the optimizer of \eqref{eq:10}
$$x^*=d\lambda_1^{-\frac12} W v_1$$
by taking the transformation $z=QW^{-1}x$. Furthermore, the positivity of $x$ can be guaranteed by the Perron-Frobenius theorem \cite{Horn85} for the non-negative and irreducible matrix $H$.

The optimal $t^*$ is obtained by the relation
$$ t^*=Tx^*/d = \lambda_1^{-\frac12} TW v_1 = TW v_1/\|TW v_1\|.$$
The proof is done.
\end{proof}

\begin{remark}
The formulation \eqref{eq:timescale} realizes the inference of model \eqref{eq:prop1} through maximum likelihood estimation. The rescaling parameter $x_g$ corresponds to the inverse splicing rate $\beta_g^{-1}$, and the normalization $\|t\|=1$ is to fix the undetermined global time scale of the whole system. The constant $d\lambda_1^{-1/2}$ in $x^*$ is not important but the orientation $Wv_1$ is essential.
\end{remark}

Our second proposal is slightly different from the first one, and it is based on the model
\begin{equation}\label{eq:prop2}
t_{cg}=t_c\beta_g+\epsilon_{cg},\quad \epsilon_{cg}\sim N(0,\sigma^2) \quad \textrm{for } c=1,\ldots,n;\ g=1,\ldots,d.
\end{equation}
With this setup, we can directly determine $\beta$ and $t$ through the following maximum likelihood formulation.

\begin{proposal}[Inference with Additive Noise]
The gene-shared latent time $t$ and the splicing rate $\beta$ can be determined by solving the minimization problem
\begin{equation}
(\beta^*,t^*) = \argmin_{\|\beta\|=1;\beta\succ 0,t\succeq0}\Vert T- t \beta^{\mathrm{T}}\Vert_F^2.
\label{eq2:timescale}
\end{equation}
\end{proposal}

\begin{thm}\label{thm:2} Assume that the inferred gene-specific cell time matrix $T$ satisfies the condition \eqref{eq:Irreduc}.
Then the optimization problem \eqref{eq2:timescale} has the unique solution $\beta^*=v_1$, where $v_1$ is the $\ell^2$-unit eigenvector corresponding to the maximal eigenvalue $\lambda_1$ of
\begin{equation}
H=T^{\mathrm{T}}T
\end{equation}
and it has positive components. The global gene-shared common time
$$t^*=Tv_1.$$
\end{thm}
\begin{proof}
Note that when $\beta$ is fixed, the optimization
$$\min_{t\in\mathbb{R}^d}\Vert T-t\beta^{\mathrm{T}}\|^2_F$$
is a least squares problem, and the minimum point is $t=T\beta/\Vert\beta\Vert^2$. Substituting it back and ignoring the normalization and positivity constraints on $\beta$ at first, we obtain
$$\min_{\beta\in \mathbb{R}^d}\Big\| T-\frac{T\beta\beta^{\mathrm{T}}}{\Vert\beta\Vert^2}\Big\|_F^2\quad \Longleftrightarrow \quad
\max_{\beta\in \mathbb{R}^d} \frac{\beta^{\mathrm{T}}T^{\mathrm{T}}T\beta}{\Vert\beta\Vert^2}$$
The rates $\beta$ can be determined up to a multiplicative constant. So we naturally take the normalization $\|\beta\|=1$ and consider the equivalent problem
\begin{equation}
\beta^* = \argmax_{\|\beta\|=1,\beta\succ 0} \beta^{\mathrm{T}}T^{\mathrm{T}}T\beta.
\label{eq:optm2}
\end{equation}

By the condition \eqref{eq:Irreduc} and the Perron-Frobenius theorem applied to the matrix $H=T^\TT T$, the optimizer of \eqref{eq:optm2} is unique and characterized by the unit eigenvector $v_1$ associated with the maximal eigenvalue $\lambda_1$ of $H$, and it has positive components.
\end{proof}

With the above proposals, we get the splicing rates $\beta^*$ and gene-shared latent time $t^*$. We can make the rescaling
$$
(\alpha_g,1,\gamma_g; t_{cg})\longrightarrow (\alpha_g \beta_g^*, \beta_g^*, \gamma_g\beta_g^*; t_{cg}/\beta^*_g),\quad g=1,\ldots,d
$$
to get more reasonable parameters with the obtained $\beta^*$.

\begin{remark}
In actual computations, when $\|t_{\bullet g}\|=0$ for some $g$, this gene will be skipped in the computation. The condition \eqref{eq:Irreduc} is not stringent if the dropout effect is not significant.  The normalization $\|t\|=1$ in \eqref{eq:timescale} is to ensure the existence of positive rates $\beta$. Another choice $\|x\|=1$ may not guarantee such positive solution.

The key difference between Proposals 1 and 2 is that they have the following comparative form
$$t_{cg}=t_c\beta_g+\beta_g\epsilon_{cg}\ (\text{Proposal 1}),\ \  t_{cg}=t_c\beta_g+\epsilon_{cg}\ (\text{Proposal 2}).$$
That is why we call Proposal 1 the multiplicative noise case, while Proposal 2 the additive noise case. It is not clear a priori which choice is more reasonable in practical situations.
\end{remark}

Both of the above two proposals assume that all of the gene expressions start from a common initial time, which is defined as 0. This assumption may be too strong since the expression for different genes may start from different instants. An extension of this point and general consideration of dropout effect for the gene-shared latent time are studied in \cite{Wang23Preprint}.

\begin{figure}[htbp]
\centering
\includegraphics{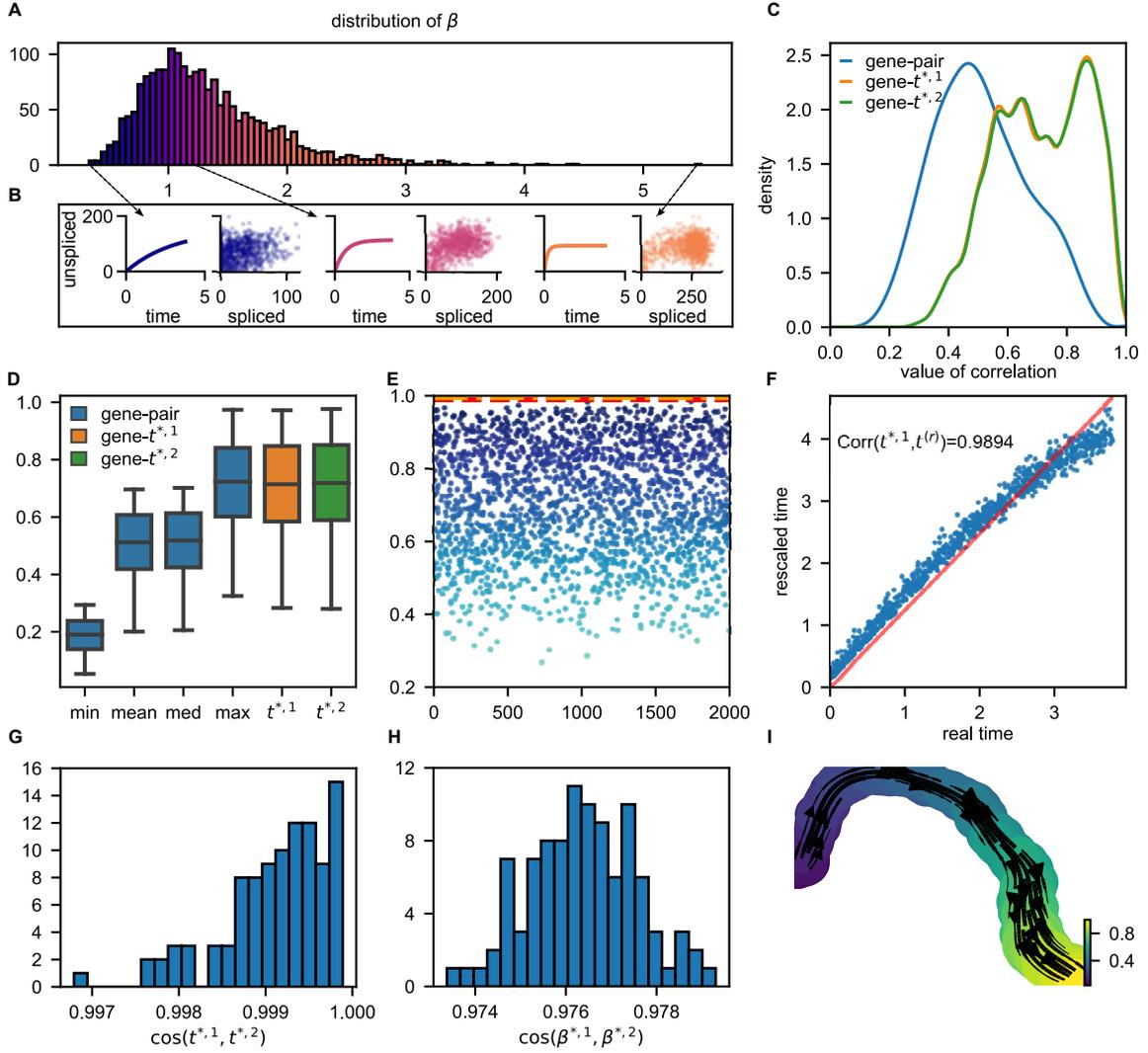}
\caption{\textbf{Resolving scaling parameters in simulation data}. (A) The empirical distribution of simulated splicing rates $\beta_g$, in which the samples are generated from a log-normal distribution. (B) The dynamics of unspliced mRNA and the simulated data with different splicing rates $\beta_g$. (C) The comparison between gene-specific time and gene-shared time. The blue curve denotes the distribution of the Pearson correlation between each pair of gene-specific time sequences, showing that the correlation between direct inferred times is low. The orange and green curves are the distribution of correlation between $t^{*,1}$, $t^{*,2}$ and all the gene-specific times, respectively, revealing a more compatible pattern. (D) A further comparison between gene-specific time and gene-shared time. The blue boxes correspond to the distributions of five important statistics of each gene-specific time sequence's correlations. The distribution of the correlation between gene-shared time and gene-specific time is drawn as the orange and green boxes. The rescaled time is as consistent as the highest correlation of genes. (E) The correlation of genes-specific time with $t^{(r)}$ compared with the correlation between the gene-shared time and $t^{(r)}$ (the red dashed line for $t^{*,1}$ and the orange one for $t^{*,2}$). Rescaled time is the most consistent time sequence to $t^{(r)}$, while some of the gene-specific times highly differ from the ground truth. (F) Scatter plot of $t^{(r)}$ and the rescaled time $t^{*,1}$, the red line denotes the fitted line. The rescaled time shows an obvious linear pattern in regard to the real time. (G) The cosine correlation between $t^{*,1}$ and $t^{*,2}$. The gene-shared times rescaled by two proposals are very close to each other. (H) The cosine correlation between $\beta^{*,1}$ and $\beta^{*,2}$, showing very similar rescaled parameters. (I) Visualization of simulated data and the fitted streamlines based on UMAP \cite{mcinnes2018umap}. The coloring is based on the rescaled time normalizing to the interval $[0,1]$, which is consistent with the embedded streamlines.}
\label{fig:rescale}
\end{figure}

\subsection{Numerical validation}

To verify the effectiveness of the proposals considered in Sec.~\ref{sec:2.2}, we make an illustration with a synthetic example.
We simulated $1000$ cells with $2000$ genes in the on stage by first sampling the parameters $(\alpha_g,\beta_g,\gamma_g)$, whose distribution is set to be ${\rm lognormal}(\mu,\Sigma)$, in which $\mu=[5,0.2,0.05]$, $\Sigma_{11}=\Sigma_{22}=\Sigma_{33}=0.16$, $\Sigma_{12}=\Sigma_{21}=0.128$, and $\Sigma_{23}=0.032$ (Fig. \ref{fig:rescale}A). This results in a typical scale of $100$ for the simulated mRNA counts. To avoid the case that the system is almost at steady state and the majority of fluctuations are caused by the observation noise, we sampled the physical real time $t^{(r)}=(t^{(r)}_c)_c$ for the cells from a uniform distribution $\mathcal{U}[0,T]$ with $T$ determined as the median of $\tau_g:=2\ln(10)/\beta_g$ for $g=1,\ldots,d$, where the number $2\ln(10)$ in $\tau_g$ is chosen such that $u(\tau_g)\approx 0.99\alpha_g/\beta_g$ which is close to the steady state. Then we computed the exact expression number by Eq.~\eqref{deter_model}, and added a Gaussian noise with mean $0$ and standard deviation $30$ to form the synthetically measured data (Fig.~\ref{fig:rescale}B). In the inference stage, We first inferred the parameters by setting the splicing rates $\beta_g=1$, then determined the time-scale parameters by the proposed methods in previous subsection. We call the optimized gene-shared common time $t^{*,1}=(t_c^{*,1})_c$ and $t^{*,2}=(t_c^{*,2})_c$ obtained from Proposals 1 and 2, respectively, and the corresponding optimized splicing rate $\beta^{*,1}$ and $\beta^{*,2}$.

As we cannot recover the physical time $t^{(r)}$of cells due to the scale invariance and an undetermined global timescale, a good rescaling method should improve the linear correlation between the inferred gene-shared time $t^*$ and the gene-specific time $(t_{\bullet g})$. This point is shown in Figs. \ref{fig:rescale}C and \ref{fig:rescale}D, in which we can find that the correlation coefficient distribution and its statistics for the correlation between $t^*$ and $(t_{\bullet g})$ for different $g$ have significant improvements compared with those for the gene pairs $(g_1,g_2)$ (i.e., the correlations ${\rm Corr}(t_{\bullet g_1},t_{\bullet g_2}) = t_{\bullet g_1}\cdot t_{\bullet g_2}/\|t_{\bullet g_1}\| \|t_{\bullet g_2}\|)$.

It is also expected that the inferred gene-shared time $t^*$ has better correlation with the real time $t^{(r)}$ than the gene-pair correlations. This is verified in Fig.~\ref{fig:rescale}E, where we can find that the $(t^{*,1},t^{(r)})$ correlation achieves a high value of $0.9894$, which is far bigger than the correlations between gene pairs. Furthermore, the scatter plot of $(t^{(r)}_c, t^*_c)$ for different cells in Fig.~\ref{fig:rescale}F shows an evident linear relation, and this linear dependence is better at early stage of the gene expression, and slightly deteriorates in later stage when the expression reaches steady states. Similar pattern can be also observed in the off stage and we omit it.

To understand the relation between the optimized time obtained from two proposals, we perform $100$ times of independent simulations by the same workflow, thus obtain $100$ pairs of $(t^{*,1}, t^{*,2})$ and $(\beta^{*,1},\beta^{*,2})$. In Figs.~\ref{fig:rescale}G and \ref{fig:rescale}H, we present the distribution of correlation coefficients for the optimized gene-shared time $t^*$ and splicing rates $\beta^*$ from two proposals, respectively. It shows that the two proposals give very close results in terms of the cosine correlations, which are around 0.999 for $t^*$ and 0.976 for $\beta^*$. This suggests both options are acceptable choices. In Fig.~\ref{fig:rescale}I, we present the streamline plot of the inferred RNA velocity with the UMAP representation and the smoothed cell coloring according to the gene-shared latent time $t^*$, which shows nice consistency between the developing flow and time progression of cells.

The issue of being unable to fix the undetermined global time-scale by the proposed approaches is due to the intrinsic drawback of the current experimental techniques. Resolution of this issue depends on further progress of the sequencing technology to extract the temporal information, such as the recent metabolic labeling technique \cite{qiu2022mapping}.

\section{Uncertainty quantification of RNA velocity}

In previous section, we proposed a method to unify the time scale between different genes which is critical to the complete parameter inference of RNA velocity models. After parameters are determined, it is also important to  evaluate the quality as well as quantify the uncertainty of the inferred parameters and computed RNA velocity. Therefore we will study the confidence interval construction of RNA velocity models through the Fisher information approach and SEM (supplemented EM) algorithm\cite{MengUsing1991}.

\subsection{Problem setup}

 For the observed data $x_{\text{obs}}=(x_{cg})_{cg}=(u_{cg},s_{cg})_{cg}$, we want to maximize the log-likelihood
 $$
 \begin{aligned}
 L(\theta| x_{\text{obs}})&=\log p(x_{\text{obs}}|\theta)\\
 &=\log\left[\prod_{cg}\int_\mathbb{R} p\left(x_{cg},t|\theta\right)\dd t\right],
 \end{aligned}
 $$
 where $\theta=(\alpha_g, \beta_g, \gamma_g)_g$, and $p(x,t|\theta)$ is the joint distribution of $(x,t)$ when $\theta$ is fixed. The marginal distribution $\int_\mathbb{R} p(x,t|\theta)\dd t$ is also called the occupancy distribution of cells in \cite{gorin2022rna}. In general $p(x,t|\theta)$ has the form
 $$p(x,t|\theta) = p(x|t, \theta)p(t|\theta)$$
 where $p(t|\theta)$ is the assumed distribution of the physical time of cells in the considered snapshot data. A working assumption on $p(t|\theta)$ is the natural choice $p(t|\theta)\equiv p(t)=\chi_{[0,T]}(t)/T$, i.e., the uniform distribution on $[0,T]$, which is independent of the parameter $\theta$.

 In the inference process,  we assume the observation noise is Gaussian with mean $0$ and variance $\sigma^2$ for all cells and genes. Then, the  log-likelihood  is
\begin{equation*}\label{eq:likelihood}
  L(\theta| x_{\text{obs}})=\log\left[\prod_{cg}\int_0^T\frac{1}{2 \pi \sigma^2} \exp \left(-\frac{\left\|x_{cg}-x_{cg}\left(t_{cg}; \theta_g\right)\right\|^2}{2 \sigma^2}\right) \cdot \frac{1}{T}\mathrm{d} t_{cg} \right].
\end{equation*}
upon taking the independent-$t$ model discussed in \cite{Li1}. From the analysis in \cite{Li1}, we know that
\begin{equation}\label{eq:CompDataDist}
\log (p(x_{cg},{t_{cg}}|\theta))= -\left\|x_{cg}-x_{cg}\left(t_{cg}; \theta_g\right)\right\|^2 +C
\end{equation}
and
$$p\left(t_{cg}|x_{cg},\theta\right)\propto \exp \left(-\frac{\left\|x_{cg}-x_{cg}\left(t_{cg}; \theta_g\right)\right\|^2}{2 \sigma^2}\right).$$
Considering $t$ as the latent variable and utilizing the EM algorithm, we have
 $$
 \theta^{(k+1)}=\underset{\theta}{\operatorname{argmin}} \int_0^T\cdots\int_0^T\sum_{cg}\|x_{cg}-x_{cg}(t_{cg} ; \theta_g)\|^2 \exp \left(-\frac{\left\|x_{cg}-x_{cg}(t_{cg} ; \theta^{(k)}_g)\right\|^2}{2 \sigma^2}\right) \prod_{cg}\mathrm{d} t_{cg}.
 $$
As $\sigma\rightarrow 0$, by Laplace asymptotics, we obtain
\begin{align}
\textrm{E-Step:} &\quad  t_{cg}^{(k)}=\underset{t}{\arg \min }\left\|x_{cg}-x_{cg}(t ; \theta_g^{(k)})\right\|^2, \label{eq:Estep}\\
\textrm{M-Step:} &\quad \theta^{(k+1)}=\underset{\theta}{\operatorname{argmin}}\sum_{cg}\left\|x_{cg}-x_{cg}(t_{cg}^{(k)} ; \theta_g)\right\|^2. \label{eq:Mstep}
\end{align}

The EM algorithm to infer the parameters is performed by iteratively estimating the rates $\theta$ and latent time $t$ through \eqref{eq:Estep} and \eqref{eq:Mstep} until convergence. Our goal in this section is to quantify the uncertainty of the inferred parameters.

\subsection{Confidence interval construction through Fisher information}

The uncertainty of the maximum likelihood estimator (MLE) can be quantified based on the classical theory of point estimation \cite{Lehmann1998}. For independent and identically distributed data,  the MLE $\hat{\theta}_n$ obtained from $n$ samples $\{x_i\}_{i=1:n}$ converges to the true parameter $\theta^{*}$ under suitable regularity conditions on $P$ in the following sense
\begin{equation}\label{eq:MLE-CLT}
\sqrt{n}(\hat{\theta}_n-\theta^{*}) \stackrel{d}{\rightarrow} \mathcal{N}\left(0, I^{-1}\left(\theta^{*}\right)\right) \text { as } n \rightarrow \infty,
\end{equation}
where
\begin{equation}\label{eq:FisherIM}
I(\theta)=-\int \nabla^2_{\theta}\log p(x|\theta) p(x|\theta)\dd x
\end{equation}
is the Fisher information matrix, and the convergence ``$\stackrel{d}{\rightarrow}$" holds in the sense of distribution. Hence, for large enough $n$, the error is approximately normally distributed
$
(\hat{\theta}_n-\theta^{*}) \stackrel{d}{\approx} \mathcal{N}\left(0,I^{-1}\left(\theta^{*}\right)/n\right),
$
which means that the $\hat{\theta}_n$  converges to $\theta^{*}$ with the error of magnitude $1/\sqrt{n}$ and a constant characterized by the inverse of the Fisher information matrix at the true value $\theta^{*}$. In practical computations, the uncertainty of the estimator $\hat{\theta}_n$ can be quantified based on approximating the Fisher information matrix $I(\theta^*)$ by its empirical form
$$\hat{I}(\theta^*|x_{\text{obs}})\approx\hat{I}(\hat{\theta}_n|x_{\text{obs}}) := -\frac1n\sum_{i=1}^n \nabla_\theta^2\log p(x_i|\hat{\theta}_n).$$

For the problem involving latent variables, the calculation of the Fisher information matrix $\hat{I}(\hat{\theta}_n|x_{\text{obs}})$ is not straightforward since the computation of the probability $p(x|\theta)$ involves the integral with respect to the latent time $t$.  Fortunately, this issue has been studied in \cite{MengUsing1991}, and the proposed approach can be utilized to approximate $I^{-1}(\theta^*)$ directly.

Following \cite{MengUsing1991}, we define the empirical information matrix $\hat{I}_{o}$ with the observed data $x_{\mathrm{obs}}$ as
$$
\hat{I}_{o}\left(\theta | x_{\mathrm{obs}}\right)=-\frac1n\nabla^2_{\theta} L(\theta| x_{\text{obs}})=-\frac1n\sum_{i=1}^n\nabla^2_\theta \log p(x_i|\theta).
$$
and its inverse at $\theta=\theta^*$ (if the inverse exists)
$$\hat{V}(\theta^*)=\left(\hat{I}_{o}\left(\theta^{*} | x_{\mathrm{obs}}\right)\right)^{-1}.$$
As shown in \eqref{eq:MLE-CLT}, the matrix $\hat{V}^*$ characterizes the uncertainty of the estimated parameter $\hat{\theta}_n$.
We can further define the complete-data information matrix with partial observable
$$
\hat{I}_{oc}(\theta | x_{\text{obs}},t)=-\frac{1}{n}\nabla^2_{\theta} L(\theta\mid x_{\text{obs}},t)
=-\frac1n\sum_{i=1}^n\nabla^2_\theta \log p(x_i,t|\theta).
$$
It is usually a simple function (e.g., Eq.~\eqref{eq:CompDataDist}), whose expectation about the conditional distribution $p(t|x_{\text{obs}},\theta)$ evaluated at $\theta=$ $\theta^{*}$ is:
$$
\hat{I}_{o c}=\left.\mathbb{E}_{t|x_{\mathrm{obs}}, \theta}\left[\hat{I}_{oc}(\theta | x_{\text{obs}},t) \right]\right|_{\theta=\theta^{*}} = -\frac1n\sum_{i=1}^n\int \nabla^2_\theta \log p(x_i,t|\theta^*) p(t|x_i,\theta^*)dt.
$$

From  \cite{MengUsing1991}, the EM algorithm (\ref{eq:Estep})-(\ref{eq:Mstep}) can be viewed as a mapping $\theta \rightarrow M(\theta)$ from the parameter space to itself, which has the form
$$
\theta^{(k+1)}=M(\theta^{(k)}), \quad \text { for } k=0,1, \ldots
$$
If $\theta^{(k)}$ converges to $\theta^*$ in the parameter space and $M(\theta)$ is continuous, then we have $\theta^*=M\left(\theta^*\right)$. By Taylor expansion in the neighborhood of $\theta^*$, we get
$$
\theta^{(k+1)}-\theta^* \approx J_M\cdot(\theta^{(k)}-\theta^*), \quad\left(J_M\right)_{i j}=\left.\left(\frac{\partial M_i(\theta)}{\partial \theta_j}\right)\right|_{\theta=\theta^*}.
$$

With the formula of total  probability $p(x_{\text{obs}},t|\theta) = p(x_{\text{obs}}|\theta)p(t|x_{\text{obs}},\theta)$, we have
\begin{equation}\label{eq:TotalProb-Decomp}
\log  p(x_{\mathrm{obs}}|\theta)=\log p(x_{\mathrm{obs}},t|\theta)-\log p(t|x_{\mathrm{obs}},\theta).
\end{equation}
Taking expectation to both sides of \eqref{eq:TotalProb-Decomp} with respect to $p(t|x_{\mathrm{obs}},\theta)$ , we get
$$
\hat{I}_o\left(\theta^*|x_{\mathrm{obs}}\right)=\hat{I}_{o c}-\hat{I}_{o m}=\hat{I}_{o c}\left(I-\hat{I}_{o c}^{-1}\hat{I}_{o m} \right)
$$
where
$$\hat{I}_{o m}:=\frac1n\left.\mathbb{E}_{t|x_{\mathrm{obs}}, \theta}\Big[-\nabla_\theta^2 \log p\left(t|x_{\mathrm{obs}}, \theta\right)\Big]\right|_{\theta=\theta^*}$$
is the missing information. The key observation in \cite{MengUsing1991} is that $J_M=\hat{I}_{o c}^{-1} \hat{I}_{o m} $, which is illustrated as below.

Implementation of the EM iterations from $\theta^{(k)}$  to $\theta^{(k+1)}$ is usually performed by taking the maximization of $Q(\tilde{\theta}|\theta):=\int L(\tilde\theta | x_{\textrm{obs}}, t) p\left(t | x_{\textrm{obs}}, \theta\right) dt$ with respect to $\tilde\theta$, i.e., we have
\begin{equation}\label{eq:g-eq}
g(\theta^{(k+1)}, \theta^{(k)}):=\int \nabla_\theta L(\theta^{(k+1)} | x_{\mathrm{obs}}, t) p(t|x_{\textrm{obs}}, \theta^{(k)}) d t=0.
\end{equation}
Generally denote \eqref{eq:g-eq} as $g(\tilde\theta,\theta)=0$ where $\tilde\theta=M(\theta)$. We can further take derivative of $g(\tilde{\theta}, \theta)$ with respect to $\theta$ to obtain
$$
\frac{\partial g(\tilde{\theta}, \theta)}{\partial \tilde{\theta}} \frac{\partial M}{\partial \theta}+\frac{\partial g(\tilde{\theta}, \theta)}{\partial \theta}=0.
$$
This leads to
\begin{equation}\label{eq:JM=IOC-IOM}
\frac{\partial M}{\partial \theta}\Big|_{\theta^*}=-\frac{\partial g(\tilde{\theta}, \theta)}{\partial \tilde{\theta}}\Big|_{(\theta^*, \theta^*)} ^{-1}\frac{\partial g(\tilde{\theta}, \theta)}{\partial \theta}\Big|_{(\theta^*, \theta^*)}.
\end{equation}
Some algebraic manipulations show that
$$
\frac{\partial g(\tilde{\theta}, \theta)}{\partial \tilde{\theta}}\Big|_{(\theta^*, \theta^*)}=-n \hat{I}_{o c} \quad \text { and }\quad \frac{\partial g(\tilde{\theta}, \theta)}{\partial \theta}\Big|_{\left(\theta^*, \theta^*\right)}=n\hat{I}_{o m}.
$$
Substitute these relations into \eqref{eq:JM=IOC-IOM}, we get $J_M=\hat{I}_{o c}^{-1} \hat{I}_{o m}$, and the uncertainty covariance matrix
$$
\hat{V}(\theta^*)= (I-J_M)^{-1} \hat{I}_{o c}^{-1}
$$

In practical computations, $\theta^*$ should be replaced with $\hat\theta_n$, i.e., the convergence value of EM iterations. The Jacobian $J_M$ at $\theta=\hat\theta_n$ can be approximated by simple difference quotient strategy with a prescribed suitable step size. In most cases, we only care about the diagonal components $\hat{v}^*_{ii}$ of $\hat{V}(\hat\theta_n)$ since they are directly related to the variances of the components $\hat\theta_{n,i}$ for $i=1,\ldots,d$. According to \eqref{eq:MLE-CLT}, we have an approximately $95\%$ confidence interval
$$
\left(\hat\theta_{n,i}-1.96 \sqrt{\frac{\hat{v}^*_{ii}}{n}},\ \hat\theta_{n,i}+1.96 \sqrt{\frac{\hat{v}^*_{ii}}{n}}\right)\quad \mbox{ for }\ i=1,\ldots,d
$$
such that $\theta^*_i$ falls in this interval. One important issue is that $\hat{V}^*$ should be positive definite according to its probabilistic meaning. However, it is not guaranteed automatically. some discussions about this point can be referred to \cite{spall2005monte,meng2017efficient}.

\subsection{Numerical validation}

Next, we validate the constructed  confidence interval for RNA velocity model using simulation data. We mainly test the accuracy of $(\alpha, \gamma)$ estimation with different parameter settings under various stages.
Steady states of active transcription and inactive silencing can be reached when the induced and repressed transcriptional phases last long enough, respectively.
However, these steady states are often difficult to capture, and most processes enter the next life process without reaching the steady state.
Here we use models that have not reached the steady state for parameter estimation.

For the on-stage, from \cite{Li1}, we know that the analytical solution of system (\ref{deter_model}) is
\begin{equation}\label{eq:on_stage1}
\begin{aligned}
&u(t)=u_0 e^{-\beta t}+\frac{\alpha}{\beta}\left(1-e^{-\beta t}\right), \\
&s(t)=s_0 e^{-\beta t}+\frac{\alpha}{\beta}\left(1-e^{-\beta t}\right)-\left(\alpha-\beta u_0\right) t e^{-\beta t}
\end{aligned}
\end{equation}
for $t<t_s$. In order to verify the sensitivity of different transcription stages to changes in splicing kinetic parameters, we simulated $n=800$ cells with $d=20$ genes in the on-stage. For convenience, we chose $\beta=1$ to avoid scale invariance issue and generate $20$ pairs of $(\alpha_g, \gamma_g)_g$ in which $\alpha=20:0.5:29.5$ and $\gamma=1.5:0.05:2.45$ to test the construction of confidence interval. The physical times of cells were sampled from a uniform distribution $\mathcal{U}[0, T]$ with $T=2\ln(10)$ to avoid the case that the system reaches and stays at the steady state.

\begin{figure}[!ht]
    \centering
    \includegraphics{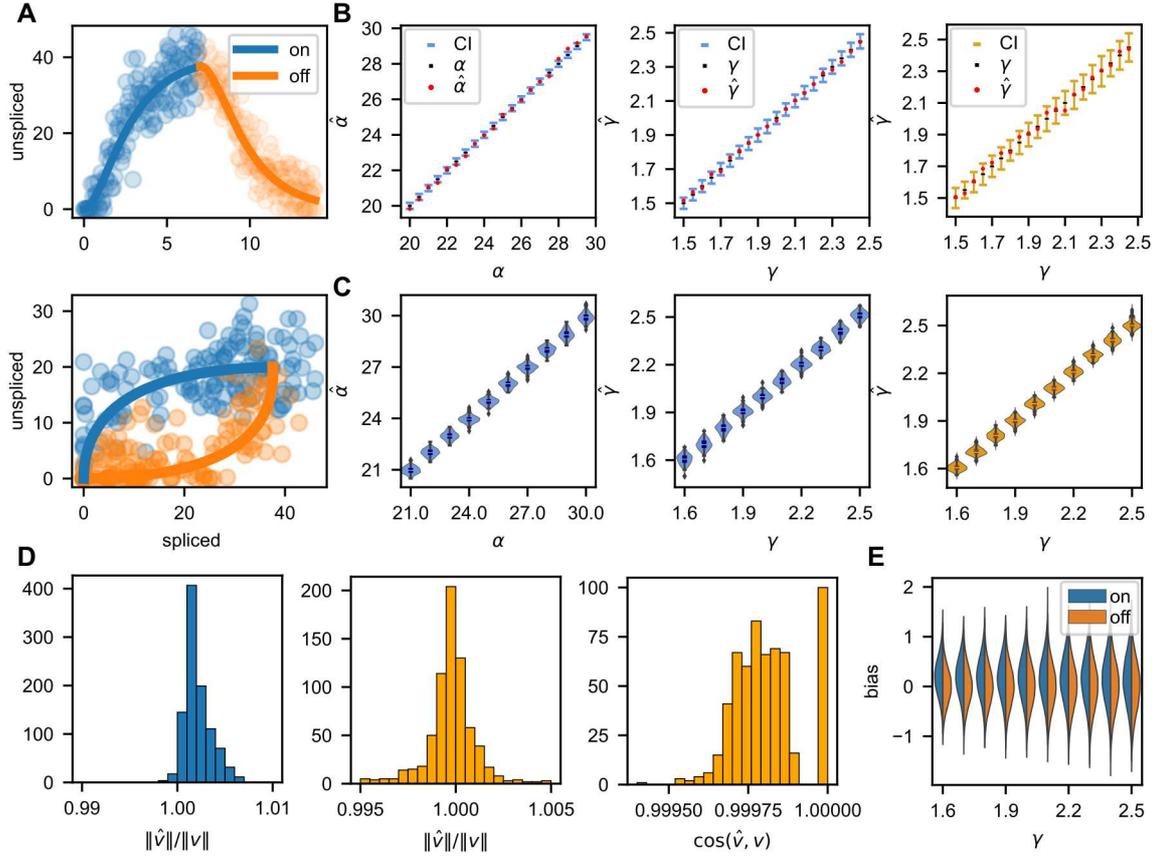}
    \caption{\textbf{Uncertainty quantification of RNA velocity using simulation data}. (A) Simulation of transcriptional process captures transcriptional induction and repression (``on" and ``off" stages) of unspliced and spliced mRNA.
   The top panel shows the abundance of spliced mRNA and the bottom panel shows unspliced and spliced mRNA in phase space. To distinguish the results, we use blue plots for on-stage and orange for off-stage.
     (B)  The  $95 \%$-prediction intervals are presented together with the inferred parameters (``on" and ``off" stages).
      It can be found that almost all of the parameters we infer are within the confidence interval.
    (C) The violin plots depict the distribution of inferred parameters of on-stage (left panel and middle panel) and off-stage (right panel). Fitted parameters mostly lie in a small range.
    (D) The norm of inferred RNA velocity is compared to the norm of true velocity at on-stage (left panel) and off-stage  (middle panel) which show high concentration near $1$, also, the cosine of the angle between inferred velocity and true velocity is shown for off-stage (right panel) which also distributed near $1$.
    (E) The bias of inferred velocity under different degradation rates $\gamma$ at on-stage and off-stage, which is close to normal distribution.
    }
    \label{fig:3}
\end{figure}

For the off-stage, from \cite{Li1}, we know that the analytical solution of off-stage is
\begin{equation}
\begin{aligned}
&u(t)=u_s e^{-\beta\left(t-t_s\right)}, \\
&s(t)=s_s e^{-\gamma\left(t-t_s\right)}-\frac{\beta u_s}{\gamma-\beta}\left(e^{-\gamma\left(t-t_s\right)}-e^{-\beta\left(t-t_s\right)}\right)
\end{aligned}
\end{equation}
for $t>t_s$. We again simulated $800$ cells with $20$ genes. The parameters were set to be the same with that in the on stage. The physical times of cells were also sampled from a uniform distribution $\mathcal{U}[0, T]$ with $T=2\ln(10)$.

The observed data was generated by adding Gaussian  noise to the dynamics, that is, $x_{\mathrm{obs}}=x_{\mathrm{true}}+\xi$ with $\xi=\mathrm{normrnd}(\mu,\sigma,2,n)$.
In the specific calculation process,  $u_{\mathrm{obs}}$ and $s_{\mathrm{obs}}$ are $n\times d$ matrices, where the observations $u_{\mathrm{obs}}^{i}$ and  $s_{\mathrm{obs}}^{i}$ are the elements representing the $i$-th column in the corresponding matrix, i.e. the unspliced and spliced mRNA for a particular gene, and the $u_{\mathrm{obs}}$ and $s_{\mathrm{obs}}$  produced have the same variance  with $\mu=0$ and $\sigma=0.2$.
Fig.~\ref{fig:3}B (left panel and middle panel) shows the inference results of parameters $\alpha$ and $\gamma$ at on-stage, respectively, together with the confidence interval of the parameters. Fig.~\ref{fig:3}B (right panel) shows the inference results of parameter $\gamma$. It can be seen from Fig.~\ref{fig:3}B that reasonable inference results can be obtained for Gaussian noise in both on-stage and off-stage.

  In order to further validate the constructed confidence interval, we show the distribution of the inference parameters $\alpha$ and $\gamma$ in Fig. \ref{fig:3}C. To better present the results, we only use $10$ genes here whose parameters are $\alpha_g = 21:1:30$, $\beta_g=1$ and $\gamma_g=1.6:0.1:2.5$. We added Gaussian noise and performed inference for each pair of parameters for $100$ times to show the violin plot.
  As can be seen from these figures, the estimated parameters closely match the true  quantities which shows the reliability of the inference procedure. Also, this result can validate the confidence interval as shown in Fig.~\ref{fig:3}B.

  In order to show the reliability of the inferred RNA velocity, we compared the norm of the inferred RNA velocity with the norm of the actual velocity, and the cosine value of the angle between the two velocities in Fig.~\ref{fig:3}D.
  We  randomly sampled 100 pairs of log-normally distributed parameters $(\alpha_g,\gamma_g)$, i.e., $\theta=(\alpha,\gamma)$ with $\log(\theta)=N(\mu_1,\Sigma)$, where $\mu_1=(3, 0.15)$ and $\Sigma=0.1I_2$. We assumed that the observation duration $T = 2$ and cell times were randomly sampled from $\mathcal{U}[0, T]$ with noise $\xi$,   i.e., we used 100 genes and 800 cells to infer the RNA velocity.  It can be seen from Fig.~\ref{fig:3}D that the ratio of the norms and the cosine value of the velocity angles are distributed around 1, which indicate that our inferred velocity size and direction are reliable.

We denote our inferred velocity as $\hat{v}$. Through the definition of RNA velocity $v^*=u-\gamma s$ when assuming $\beta_g=1$ for all genes,  we have  $\hat{v}\approx u+\xi_1-\gamma (s+\xi_2)=v^*+(\xi_1-\gamma\xi_2)$, which implies that the actual and inferred velocity are not only affected by noise but also related to the selection of $\gamma$, thus the ratio of the norm of velocities is affected by $\gamma$.
 In order to demonstrate this statement, we further studied the impact of $\gamma$ and displayed the results in Fig.~\ref{fig:3}E. Here we simulated $1600$ cells with $1000$ genes in which half the cells were in the on-stage while others were in the off-stage. As we aim to test the influence of $\gamma$, we sampled $100$ values of $\alpha$ from a log-normally distribution with $\mu=3$ and $\sigma=0.1$, and constructed $1000$ genes with $\gamma=1.6:0.1:2.5$ paired with sampled $\alpha$. Then the inference method was applied to the simulated data and we obtained the bias between the true velocity and the inferred velocity. Note that we computed the bias gene-wise here, i.e., for each $\gamma$, we tested the bias of inferred velocity under different stages, various cell physical times and values of $\alpha$. It can be seen from the figure that the error of each selected component is close to normal distribution and the variance increases as $\gamma$ increases.

\section{Optimal choice of kernel bandwidth in random walk construction}

Assuming that the RNA velocity obtained accurately describes the actual dynamics locally, a key step in the downstream analysis is the construction of a cellular random walk, i.e., the Markov model on data, which reflects more global and long-time information about cell state-transition dynamics\cite{setty2019characterization,zhou2021dissecting,Lange2022}, and is also widely used to visualize the streamlines of RNA velocity and the embedding of cells in current practice\cite{MannoRNA2018,BergenGene2020,atta2022veloviz}. In our previous work \cite{Li1}, we derived the continuum limits (i.e., differential equations) of cellular random walk induced by various RNA velocity kernels. Empirically, when dealing with single-cell data of finite or sometimes limited sample size, the hyper-parameters (especially the bandwidth $\epsilon$ in Gaussian kernel or number of neighbors $k$ in kNN kernel) in random walk construction have a significant impact on the quantitative behavior of dynamics and downstream analysis. We will study the effects of  hyper-parameters and sample size on the random walk convergence rate, elucidating their optimal choice in practice and gaining insights for algorithm implementation.  For simplicity, we will only focus on the choice of the optimal kernel bandwidth $\epsilon$ for Gaussian-cosine scheme. The other cases can be analyzed similarly.

\subsection{Problem setup}
Let $(u_i,s_i)\in \mathbb{R}^{2d}$ be the unspliced and spliced gene expression vectors of cell $i$ for $i=1,2,\ldots,n$. To define the probability of transition dynamics between different cells, the randomness introduced by extrinsic or intrinsic noise \cite{eling2019challenges,hilfinger2011separating,zhou2021stochasticity} as well as directed state-transition in relation to RNA velocity needs consideration \cite{MannoRNA2018,BergenGene2020}. The transition between two cells usually involves both drift and diffusion effects. For the diffusion part, the popular Gaussian diffusion kernel has the form
$$
d_\epsilon\left(s_i, s_j\right)=h\Big(\frac{\left\|s_i-s_j\right\|^2}{\epsilon}\Big),
$$
where the function $h(x)$ is usually chosen as a smooth function with exponential decay. For the drift part, we consider the velocity kernel $v\left({s}_i, {s}_j\right)=g\left(\cos \left\langle{\delta}_{i j}, {v}_i\right\rangle\right)$, where $\delta_{i j}=s_j-s_i$, $v_i=\beta \circ u_i-\gamma \circ s_i$ is the RNA velocity, $\langle\delta_{i j}, {v}_i\rangle$ represents the angle between $\delta_{i j}$ and ${v}_i$, and $g(\cdot)$ is a bounded, positive, and non-decreasing function. The
overall transition kernel is then defined by
$$k_\epsilon\left(s_i, s_j\right)=d_\epsilon\left(s_i, s_j\right) \cdot v\left(s_i, s_j\right).$$
And the transition probability matrix $P_{\epsilon}=\left(p_{i j}\right)_{i, j=1: n}$ among cells through the Gaussian-cosine scheme is defined by
$$
p_{i j}=\frac{k_{\epsilon}\left({s}_{i}, {s}_{j}\right)}{\sum_{j=1}^{n} k_{\epsilon}\left({s}_{i}, {s}_{j}\right)},~s_j\sim q(y),
$$
where  $\sum_{j=1}^{n} k_{\epsilon}\left({s}_{i}, {s}_{j}\right) $ are row normalization factors.

The study of the continuum operator limit of $P_{\epsilon}$ when the number of samples is assumed as infinity has been investigated in \cite{Li1} by considering the operator $\mathcal{G}_{\epsilon}$ acting on a smooth function $f$ defined as
$$
\mathcal{G}_{\epsilon} f(x)=\frac{1}{\epsilon^{\frac{d}{2}}} \int_{\mathbb{R}^{d}} k_{\epsilon}(x, y) f(y) \mathrm{d} y.
$$
From Lemma \ref{lem1} in Appendix \ref{appendixA} (i.e., Lemma 4.1 in \cite{Li1}), the operator $\mathcal{G}_{\epsilon}$ for Gaussian-cosine scheme has the expansion
\begin{equation}\label{muf}
\begin{aligned}
\mathcal{G}_{\epsilon} f(x) &=\frac{1}{\epsilon^{\frac{d}{2}}} \int k_{\epsilon}(x, y) f(y) \mathrm{d} y \\
&=m_{0} f(x)+\sqrt{\epsilon} m_{1} \hat{v}(x) \cdot \nabla f(x)+O({\epsilon}),
\end{aligned}
\end{equation}
where $m_0,m_1$ are constants depending on functions $g,h$ in the diffusion and velocity kernels (see detailed connections in Appendix \ref{appendixA}), and  $\hat{v}(x):=v(x)/\|v(x)\|$ where $v(x)$ is the RNA velocity in the continuum formulation.

Then given the sample probability density $q(\cdot)$, the continuous transition kernel has the form
$$
p_{\epsilon}(x, y)=\frac{k_{\epsilon}(x, y)q(y)}{d_{\epsilon}(x)}, \quad d_{\epsilon}(x)=\int k_{\epsilon}(x, y) q(y) \mathrm{d} y .
$$
Define the operator
$$
\mathcal{P}_{\epsilon} f(x)=\int p_{\epsilon}(x, y) f(y) \mathrm{d} y
$$
and the discrete generator
\begin{equation}
\mathcal{L}_{\epsilon}=\frac{\mathcal{P}_{\epsilon}-I}{\sqrt{\epsilon}}.
\end{equation}
From  Theorem 4.1 in \cite{Li1}, we have the convergence of the generator for the Gaussian-cosine scheme
\begin{equation}
\lim _{\epsilon \rightarrow 0+} \mathcal{L}_{\epsilon} f=\mathcal{L} f:=\frac{m_{1}}{m_{0}} \hat{v}(x) \cdot \nabla f(x), \quad \hat{v}(x):=\frac{v(x)}{\|v(x)\|}.
\end{equation}
Indeed, we can further identify the higher order expansion of $\mathcal{L}_{\epsilon}$ as
\begin{equation}\label{eq:Leps}
\mathcal{L}_{\epsilon}f(x) =\mathcal{L}f(x)+O(\sqrt{\epsilon})
\end{equation}
since
$$\begin{aligned} \mathcal{P}_\epsilon f(x)
& =\frac{\mathcal{G}_\epsilon(f q)(x)}{\mathcal{G}_\epsilon q(x)} =\frac{m_0 f(x) q(x)+\sqrt{\epsilon} m_1 \hat{{v}}(x) \cdot \nabla(f q)(x)+O(\epsilon)}{m_0 q(x)+\sqrt{\epsilon} m_1 \hat{{v}}(x) \cdot \nabla q(x)+O({\epsilon})} \\
& =f(x)+\sqrt{\epsilon} \mathcal{L}f(x)+O(\epsilon) .\end{aligned}$$

In practical computation, we merely get a limited amount of samples. Beyond investigating the convergence result of the discrete operator to the continuous infinitesimal generator in the limit $\epsilon \rightarrow 0$ when the sample size is assumed as infinity, it is also important to understand the optimal kernel bandwidth $\epsilon$ when the sample size $n$ is finite. This can be achieved by analyzing the bias and variance tradeoff of discrete models in approximation to their continuum limit.

\subsection{Estimation of operator convergence and algorithmic insight}

When the sample size $n$ is finite, the discrete generator $\mathcal{L}_{\epsilon,n}$ acting on a smooth function $f$ is defined as
\begin{equation}
\mathcal{L}_{\epsilon,n}f(x)=\frac{1}{\sqrt{\epsilon}}\Big(P_{\epsilon,n}f(x)-f(x)\Big)=\frac{1}{\sqrt{\epsilon}}\left(\frac{\frac{1}{n}\sum_{j=1}^{n}k_\epsilon(x,s_j)f(s_j)}{\frac{1}{n}\sum_{j=1}^{n}k_\epsilon(x,s_j)}-f(x)\right).
\end{equation}
 Then, we have the following estimate.
\begin{thm}[Finite sample approximation of the operator $\mathcal{L}_{\epsilon}$]\label{thm:discreteerror}
Let $s_1,s_2,\cdots,s_n$ be $n$ independent
and identically distributed samples in $\mathbb{R}^{d}$ with probability density $q(x)$. Suppose that $f\in C_0^\infty(\mathbb{R}^d)$, which is a smooth function with compact support. Then, we have the error estimate
$$
|\mathcal{L}_{\epsilon,n}f(x)-\mathcal{L}_{\epsilon}f(x)|=O\left(\frac{1}{\sqrt{n}\epsilon^{\frac{d}{4}}}\right)
$$
in the sense that both the probability
$$
p(n,\alpha):=\mathbb{P}\Big(|\mathcal{L}_{\epsilon,n}f(x)-\mathcal{L}_{\epsilon}f(x)|>\alpha\Big)
$$
and $1-p(n,\alpha)$ have the $O(1)$ magnitude in $(0,1)$ only when $\alpha=O\big(1/(\sqrt{n}\epsilon^{\frac{d}{4}})\big)$.
\end{thm}

 The proof of Theorem \ref{thm:discreteerror} will be deferred to  Appendix \ref{appendixB}. Based on Theorem \ref{thm:discreteerror} and Eq. \eqref{eq:Leps}, we obtain the estimate
\begin{align}\label{eq:bias}
    |\mathcal{L}_{\epsilon,n}f(x)-\mathcal{L}f(x)|
   &=
    | \mathcal{L}_{\epsilon,n}f(x)-\mathcal{L}_{\epsilon}f(x)+\mathcal{L}_{\epsilon}f(x)-\mathcal{L} f(x)|\nonumber \\
    &\leq  | \mathcal{L}_{\epsilon}f(x)-\mathcal{L}_{\epsilon,n}f(x)|+ | \mathcal{L}_{\epsilon}f(x)-\mathcal{L}f(x)|\nonumber \\
  &=O\left(\frac{1}{\sqrt{n}\epsilon^{\frac{d}{4}}}+\sqrt{\epsilon}\right).
\end{align}
Compared with the continuum limit result in \cite{Li1}, the $O\big({1}/({\sqrt{n}\epsilon^{\frac{d}{4}}})\big)$ term quantifies the influence of finite sample size in cellular random walk.

The above estimation suggests that to achieve the optimal approximation of $\mathcal{L}$ by $\mathcal{L}_{\epsilon,n}$, the best choice of $\epsilon$ is
\begin{equation}
\epsilon = n^{-2/(d+2)}
\end{equation}
when the sample size $n$ is fixed. This is obtained when the two error terms in \eqref{eq:bias}, the `variance' and  bias, are balanced. In this case, the optimal error of the operator approximation is $O(n^{-1/(d+2)})$.

The result also shows the ``curse of dimensionality" when using the velocity-induced cellular random walk to approximate the dynamics of continuous ODEs.
As the dimensions of the input data $d$ increases, the overall error $O(n^{-1/(d+2)})$ deteriorates, which suggests insufficient model accuracy.
This analysis leads to several possible interpretations for the downstream RNA velocity analysis:
The first insight is that to achieve a more stable inference on the cellular development path, we should better use a relatively low dimensional space for genes instead of a very high dimensional space.
The second possible interpretation is that although the operator approximation is bad in a very high dimensional space, the overall direction of developments, such as the streamline visualization of the cells along the main backbone, can still be estimated in a relatively accurate manner. The  biological pathways and highly correlated gene functional modules underlying real scRNA-seq data \cite{dong2022integrating} could also further reduce the effective dimension of the data manifold, and lead to larger convergence rate empirically than theoretical results \cite{spigler2020asymptotic}. However, this point is not straightforward to be analyzed, and we leave it to be studied in the future.

\subsection{Numerical validation}

\begin{figure}[!h]
\includegraphics{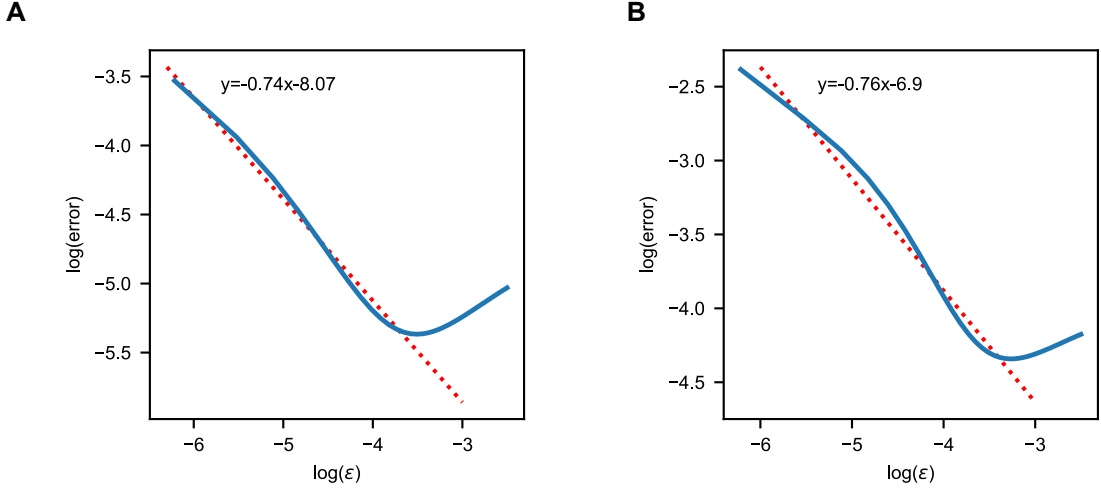}
\caption{\textbf{Effect of kernel bandwidth $\epsilon$ on operator approximation}. The logarithmic plots of the operator approximation error for the linear case $ f_1(x)=x_1+x_2+x_3$ (left panel) and nonlinear case  $f_2(x)=x_3^2$ (right panel). The errors of both cases are averaged over $n=2000 $ samples. The minimal errors are obtained when $\ln(\epsilon)\approx -3.45$  (left panel) and $\ln(\epsilon)\approx -3.20$  (right panel). When $\epsilon\lesssim n^{-2/5}$, the dominant term of the error should be $O({1}/(\sqrt{n}\epsilon^{\frac{3}{4}}))$. This gives the slope $-3/4$ in theory, which is close to the estimated value $-0.74$ (left panel) or $-0.76$ (right panel) by linear regression.}
\label{fig_error}
\end{figure}

In this subsection, we present a toy example to show the convergence rate of $\mathcal{L}_{\epsilon,n}f$. We took $d=3$ and chose a linear function $f_1(x)=x_1+x_2+x_3$ and a nonlinear function $f_2(x)=x_3^2$ to perform the numerical simulations. We first generated $n=2000$ samples $(u^{(k)},s^{(k)})_{k=1:n}=(u(t_k),s(t_k))_{k=1:n}$ according to the RNA velocity dynamics
$$\frac{du_g}{dt}=\alpha_g - \beta_g u_g,\quad \frac{ds_g}{dt}= \beta_g u_g(t) -\gamma_g s_g(t),\quad (u_g,s_g)|_{t=0}=(0,0),\quad {\rm for}\ g=1,\cdots,d$$
by choosing $t_k\sim \mathcal{U}[0,T]$ for $k=1.\cdots,n$ where $T=2\ln 10$.  Here we chose the parameters $\alpha=(20,20.5,21)^\mathrm{T}$, $\beta=(1,1,1)^\mathrm{T}$, and $\gamma=(1.5,1.55,1.6)^T$, where each component corresponds to the index $g=1,2,3$, respectively. We then generated $n=2000$ samples $\{x_k\}$ with velocity $\{v_k\}$ by setting
$$x_k=s^{(k)}+\epsilon_k,\quad \epsilon_k\sim N(0, 0.5I_3), $$
and $v_k=\beta\circ u^{(k)}-\gamma \circ x_k$ for $k=1,\cdots,n$. In the downstream analysis, we chose $g(x)=\exp(x)$, $h(x)=\exp(-x)$ and defined the root-mean-squared error as
$$
\mathrm{error}=\Bigg[\frac{1}{n}\sum_{k=1}^{n}\Big(\mathcal{L}_{\epsilon,n}f(x_k)-\mathcal{L} f(x_k)\Big)^2\Bigg]^\frac12
$$
by averaging over $n=2000$ samples. In the simulation, we chose $\epsilon=0.002:0.002:0.082$ and the results are shown in Fig. \ref{fig_error}.

According to the estimate \eqref{eq:bias}, we know that when $\epsilon\lesssim O(n^{-2/5})$, the `variance' term is dominant with the order
${1}/(\sqrt{n}\epsilon^{{3}/{4}})$. Thus, as $\epsilon$ increases, the $\ln({\rm error})$ versus $\ln(\epsilon)$ plot should present a linear relation with slope $-3/4$ theoretically. This is verified in Fig. \ref{fig_error}A and B, in which the linear fitting gives the slope $-0.74$ for the linear case shown in the left panel and $-0.76$ for the nonlinear case shown in the right panel. When $\epsilon\gtrsim O(n^{-2/5})$, the bias term is dominant and the error curve demonstrates a turn-over at $\ln(\epsilon) \sim \ln(n^{-2/5})\approx -3.04$ theoretically, which is close to the computed minimum point at $\ln(\epsilon)\approx -3.45$  in linear case and $\ln(\epsilon)\approx -3.20$  in nonlinear case.

\section{Transition time estimation among cell states}

After constructing the Markov chain among individual cells induced by RNA velocity, the downstream analysis could be performed to reflect the long-term and global (i.e., among multiple cell states) dynamics of cell-state transitions\cite{Li1,qiu2022mapping,Lange2022}. In our previous work\cite{Li1,zhou2021dissecting}, we proposed the approach based on transition path theory to infer coarse-grained cell lineage and quantify corresponding likelihoods from the cellular random walk. Another practical task is to quantify the duration of the transition from one cell to another along the transition paths, i.e. defining the pseudo-temporal distance\cite{trapnell2014pseudo} between cells induced by RNA velocity, which could be realized by the first hitting time analysis described below.

\subsection{Problem setup}

Suppose that the RNA velocity induces the cellular random walk with transition probability matrix $P=(p_{ij})$ which was defined in Section 4.1. Our goal is to define a pseudo-temporal distance $T_i^A$ from cell $i$ to the cell set $A$ which reflects the state-transition time based on the Markov chain model. When the cell set $A=\{j\}$, $T_i^A$ gives the pseudo evolution time from cell $i$ to cell $j$.

For Markov chain $\{X_n,n\geq 0\}$, the first hitting time of a set $A$ is defined as
$$
\tau^A = \inf\{n\geq0 : X_n\in A\},
$$
where $A$ is a subset of the state space. The mean first hitting time for the process to reach $A$ starting from $i$ is given by
$$
k_i^A=\mathbb{E}_i(\tau^A)=\sum_{n<\infty}n\mathbb{P}_i(\tau^A=n)+\infty\mathbb{P}_i(\tau^A=\infty)
$$
where $\mathbb{E}_i$ and $\mathbb{P}_i$ denotes the expectation and probability conditioned on $X_{0}=i$, respectively. The quantity $k_i^A$ serves as the rational proposal for pseudo-temporal distance $T_i^A$, and we will demonstrate the equations to calculate it below.

We will consider two types of hitting times to describe the evolution time from one cell to another cell set $A$. Firstly, we show that it is straightforward to use the Eq. \eqref{hitting} below to compute the mean first hitting time when there is no bifurcation. Secondly, we use a simplified model to demonstrate the limitation of this approach when there is the "bottleneck" state and cell-state differentiation in dataset. We then show how to get a biologically meaningful time by utilizing the taboo set concept.

\subsection{Transition time estimation through first hitting time analysis}

The computation of $k_i^A$ is based on the following Lemma (see, e.g., \cite[Theorem 1.3.5]{norris1998markov}).

\begin{lem}\label{lemma3}
The vector of mean first hitting times $k^A=(k^A_i)_i$ is
the minimal non-negative solution to the system of linear equations
\begin{equation}
\label{hitting}
    \left\{
    \begin{split}
    & k_i^A=0,\quad & i\in A,\\
    & k_i^A=1+\sum_{j\not\in A}p_{ij}k_j^A, \quad &i\not\in A.
    \end{split}\right.
\end{equation}
\end{lem}

To solve Eq. \eqref{hitting}, we use the following iterations
\begin{equation}
\label{hitting1}
    \begin{split}
    & K_n^A=\bm{1}+QK_{n-1}^A,\quad K_0^A=\bm{1},
    \end{split}
\end{equation}
where $K_n^A$ is the $n$th iteration of $k^A$, and $Q=(q_{ij}):=(p_{ij})_{i,j\in A^c}$, i.e., the matrix formed by removing the row and column elements corresponding to $i\in A$ from the transition probability matrix $P$. Next we show that the iteration (\ref{hitting1}) is a contraction mapping, i.e., the spectral radius $\rho(Q)<1$.

\begin{thm}
Assume the velocity-induced cellular random walk with transition probability matrix $P$ is  irreducible, then the iteration \eqref{hitting1} is a contraction mapping, i.e., $\rho(Q)<1$.
\end{thm}
\begin{proof}
Without loss of generality, assume
$$
P=\left[\begin{array}{cc}
Q & R_1 \\
R_2& B
\end{array}\right],
$$
where $Q$ and $B$ describe the transition probabilities among the states in $A^c$ and $A$, respectively. $R_1$ and $R_2$ describe the transition probabilities from the states in $A^c$ to $A$ and $A$ to $A^c$, respectively, which should not be zero since  $P$ is irreducible.

We will first consider the case that $Q$ is irreducible. In this case, if all of the row sums of the matrix $Q$ are strictly less than 1, then the conclusion holds simply by Gershgorin circle theorem \cite{Horn85}. Otherwise, we have $\rho(Q)\leq 1$ and there is at least one row of $Q$ such that its sum is strictly less than 1. Below we show that the assumption $\rho(Q)=1$ will lead to contradiction.

Utilizing the Perron-Frobenius theorem, we have the Perron vector $x$ with positive components such that
$$Qx=\rho(Q)x=x.$$
Suppose $x_l=\max\{x_k\}_{k\in A^c}$ and define $y=x/x_l$. We have $Qy=y$ and
$$\sum_{k}q_{lk}y_k=y_l=1.$$
Since $y_k\leq 1$ and  $\sum_kq_{lk}\le 1$, the above  identity requires that $y_k=1$ for $q_{lk}> 0$, i.e., the neighborhood states of $l$. We can apply similar arguments to these states $k$, which eventually lead to $y\equiv 1$. While this contracts with the condition that at least one row rum of $Q$ is strictly less than 1. So, we have $\rho(Q)<1$ when $Q$ is irreducible.

In general cases, we can decompose $A^c$ into several irreducible components and transient states. For each irreducible component, it has the transition probability sub-matrix with similar structure as the above case. Thus the spectral radius is strictly less than 1. For the transient states, the transition probability sub-matrix has the property that all row sums are less than 1, thus the spectral radius is also strictly less than 1. Overall, we have $\rho(Q)<1$ in the general cases.
\end{proof}
\begin{remark}
The result $\rho(Q)<1$ ensures that $N=(I-Q)^{-1}$ is well-defined, which is called the fundamental matrix in the literature \cite{Kemeny76}. Here we give a self-contained linear algebra proof instead of probabilistic arguments. The above theorem also tells us that $k^A=N\cdot \bm{1}$ although it is not a feasible approach to compute $k^A$ due to the ill-conditioning of $I-Q$.
\end{remark}

\begin{figure}[htbp!]
\includegraphics{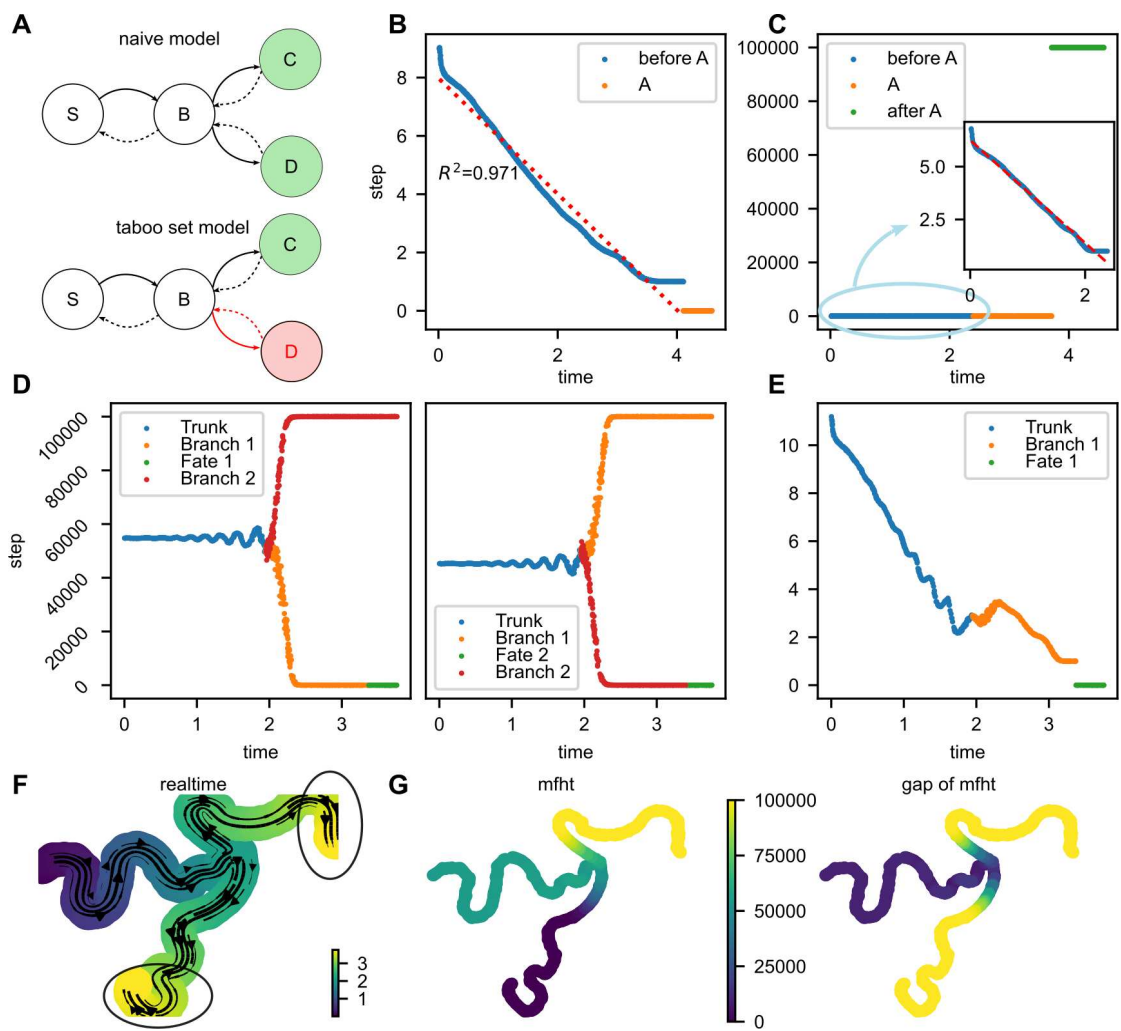}
\caption{\textbf{Estimating the pseudo-temporal distance via the first hitting times}. (A) Schematics of the naive  and taboo set models. (B) The physical time of the synthetic model compared with the mean first hitting time to the target set in the end. There is a linear pattern for the computed mean first hitting time and the red line is obtained by linear regression. (C) The physical time compared with the mean first hitting time to a target set in the middle. For cells before this target set $A$, there is still a linear pattern which is shown in the inset, in which the red dashed line is obtained by linear regression. (D) In the bifurcation case, the physical time compared with the mean first hitting time to one branch computed by iterative method. Cells in the other branch and before bifurcation have very large mean first hitting times. (E) By setting the other branch as the taboo set, the mean first hitting time shows a nearly linear pattern. (F) Streamline embedded UMAP plot of the synthetic bifurcation data. The two expected termination sets are circled. (G) UMAP of the synthetic data. The left panel is colored according to the mean first hitting time to the lower branch termination cells, and the right panel is colored according to the absolute value of the difference between two computed mean first hitting times to two circled branches in Fig.~\ref{fig_transition}F.}
\label{fig_transition}
\end{figure}

The above approach is only useful when there are no bifurcations, i.e., no differentiations in cell development. This can be illustrated by a simplified model shown in Fig. \ref{fig_transition}A, in which we model the stem cell as state $S$, the developmental bottleneck as state $B$, and two differentiated states as $C$ and $D$, denoting different fates of cell differentiation. The state $S$ can only transit to $B$, while $B$ is able to transit to $C$ and $D$, or back to $S$. We assume the transitions have some preferred directionality, i.e., it is nearly impossible to transit back along the directed developmental pathway.  The above assumption amounts to set that $p_{SB}$, $p_{BC}$, $p_{BD}$ are $O(1)$ and $p_{BS}, p_{CB}$, $p_{DB}$ are  $O(\epsilon)$. Heuristically, we can set the transition probability matrix $P$ as
$$
P = \left(
\begin{matrix}
    0 & 1 & 0 & 0\\
    \epsilon & 0 & p & q \\
    0 & \epsilon & 1-\epsilon & 0\\
    0 &\epsilon & 0 & 1-\epsilon
\end{matrix}
\right),
$$
in which $p$ and $q$ are probabilities of $O(1)$  and $p+q=1-\epsilon$.

In this setup, the mean first hitting time of each state to the target set $C$ can be obtained according to \eqref{hitting} as
\begin{equation}
    \left\{\begin{aligned}
        & k_S^C = 1+k_B^C\\
        & k_B^C = 1+qk^C_D+\epsilon k_S^C\\
        & k_D^C = 1+\epsilon k_B^C+(1-\epsilon)k_D^C,
    \end{aligned}\right.
\end{equation}
from which we get
$$k_S^C = 1+\frac{q+\epsilon+\epsilon^2}{\epsilon p},\ k_B^C = \frac{q+\epsilon+\epsilon^2}{\epsilon p},\  k_D^C = \frac{1+\epsilon^2}{\epsilon p}.
$$
We have $k^C_{S}, k^C_{B}, k^C_D\sim O(\epsilon^{-1})$, i.e., we need a prominent long transition time to reach $C$ from any other state, including $S$ and $B$, which looks counter-intuitive. The reason is that although the states $S,B$ are in the upstream of $C$ in the cell development, they have $O(1)$ probability to reach $D$ in an $O(1)$ timescale, while it is very difficult to transit back from $D$ to $B$ once $D$ is reached. This effect finally makes the overall transition time from any other states to $C$ are extremely large. However, this does not reflect the biological intuition that the transition time from $S$ and $B$ to a specific differentiated state $C$ or $D$ is in $O(1)$ timescale.

Biologically, we are mainly interested in the transition time that the stem/bottleneck cell differentiates to a specific cell type instead of letting it transit to another differentiated state and then get back. In this case, the other differentiated cell states beyond the interested states  should be  neglected and form a forbidden set. We call this set a taboo set, and compute the mean first hitting time of state $S$ or $B$ to $C$ conditional on not reaching the taboo set $H=\{D\}$. This is illustrated in the taboo set model in Fig. \ref{fig_transition}A. Denote the first hitting time with taboo set $H$ by
$$_H\tau^A = \inf\{n\geq0 : X_n\in A\ {\rm and}\ X_m\not\in H\ {\rm for}\ m\le n\}$$
and the mean first hitting time by $_H k^A$,
$$
_Hk_i^A = \sum_{n<\infty}n\mathbb{P}_i(_H \tau^A=n)+\infty\mathbb{P}_i(_H \tau^A=\infty).
$$
Then $_Hk^A$ satisfies
\begin{equation}
_Hk_i^A=1+\sum_{j\not\in A\cup H}p_{ij}{}_Hk_j^A,\quad i\not\in A\cup H.
\end{equation}
So we obtain
$$
\left\{
\begin{aligned}
    & {}_Hk_S^C = 1+{}_Hk_B^C\\
    & {}_Hk_B^C = 1+\epsilon{}_Hk_S^C,\\
\end{aligned}
\right.
$$
from which we get
$$_Hk_S^C = \frac{2}{1-\epsilon}\approx 2,\quad_Hk_B^C = \frac{1+\epsilon}{1-\epsilon}\approx 1.$$
This result reflects the intuition that the transition time from $S$ to $C$ and from $B$ to $C$ are about 2 and 1, respectively, by simply counting the transition steps in the Markov chain.  By setting the taboo set, the behavior of the tabooed process  is similar to the case when there is no bifurcation, and the taboo set model acts as a pruning strategy.

\subsection{Numerical validation}

To verify the applicability of our proposal on the evolution time estimation, for the non-bifurcation situation, we simulated $1000$ cells with $2000$ genes at on-stage to generate the synthetic data. The splicing rate $\beta_g$ was fixed to $1$ to avoid considering the scale invariance issue, and the transcription rates $\alpha_g$ and degradation rates $\gamma_g$ were sampled from a log-normal distribution with mean $\mu=[5,0.05]$ and covariance matrix $\Sigma=0.16I_2$. The physical time of cells were sampled from $\mathcal{U}[0,T]$ with $T=2\ln(10)$. After inferring the parameters, a Gaussian-cosine kernel was constructed. We chose the set $A$ as the top $100$ cells having the largest sampled real time, and compute the mean first hitting time from any cell to this set by solving \eqref{hitting}. From Fig. \ref{fig_transition}B we can find that the computed mean first hitting time matches well with the real time of cells upon ignoring a scaling constant and the R-squared of linear regression is $0.968$, which validates our proposal in the considered simple synthetic example.

To compute the mean first hitting time of non-terminal states, we applied the iteration \eqref{hitting1}, from which we can tell the relevant order of cells. We chose the cells with the rank order $500\sim 800$ as set $A$, by which the other cells were separated as two groups with weak links. By iterating for a sufficiently long time, we can tell from the result the relevant order of cells to set $A$. Shown in Fig. \ref{fig_transition}C is the iterated mean first hitting time after $1\times 10^5$ iterations, as cells after set $A$ are nearly impossible to transit back to cells prior to set $A$, the mean first hitting time of cells posterior to $A$ is very large and is the scale of iteration time. Mean first hitting times of cells prior to set $A$ still show a linear decay pattern. However, for the cells that are close to $A$, the transition times show fluctuations and an increasing trend since these cells have small probability to perform transitions to the cells posterior to $A$ due to the weak link between them. This phenomenon also partially suggests the adoption of the taboo set model in this simple case.

To test the taboo set model, we first applied the iterative method directly to a synthetic bifurcation data. Here to produce bifurcation, we sampled the parameters $(\alpha_g,\beta_g,\gamma_g)$, whose distribution was set to be ${\rm lognormal}(\mu,\Sigma)$, in which $\mu=[5,0.2,0.05]$, $\Sigma_{11}=\Sigma_{22}=\Sigma_{33}=0.16$, $\Sigma_{12}=\Sigma_{21}=0.128$, and $\Sigma_{23}=0.032$ which is the same as the setup in Section 2. Then we used the true parameters in the computation of the mean first hitting time.  The physical time of cells were sampled from $\mathcal{U}[0,T]$ with $T$ determined as the median of $\tau_g:=2\ln(10)/\beta_g$ for $g=1,\ldots,d$. The bifurcation was produced by adding the switch of gene expression from the on-stage to off-stage. For the first branch, we assigned $70\%$ of the genes to switch to off stage at $2\ln(10)/\beta_g$ and for the second branch, the rest $30\%$ genes are assigned. From Fig. \ref{fig_transition}D we can find that when setting the target set as the termination part of one branch, the mean first hitting time from another branch rapidly grows to a huge number, and the cells before the bifurcation point also have long mean first hitting times, which is similar to our analysis of the simplified naive model in previous subsection. This result implicitly indicates that we can use the naive model to detect where the bifurcation happens. As shown in Fig. \ref{fig_transition}E, by setting the second branch to be the taboo set, the mean first hitting time shows a nearly linear pattern both before and after the bifurcation point, showing the taboo set model can successfully give an estimation of transition time without considering the other bifurcation branches.

The results above show that we can use the mean first hitting time as an estimation of the physical time of cells. Such an approach for pseudo-time can be extended to real-world data in which the ground truth (i.e. physical time, bifurcation lineage, and especially taboo set) is not available. Here we will present a brief idea of the implementation in scRNA-seq data and leave further analysis and software development for future work \cite{Wang23Preprint}.

From low-dimensional visualization and velocity streamline embedding, we can identify the main fates (or branches) of differentiation as shown in Fig. \ref{fig_transition}F using existing lineage inference methods. To identify the taboo set according to the computational results, we first apply the iterative method to calculate the mean first hitting time to the expected fates on different branches and then take a postprocessing step. As shown in Fig. \ref{fig_transition}G, the left panel is colored according to the mean first hitting time to the lower branch termination cells, and the right panel is colored according to the absolute value of the difference between two computed mean first hitting times to two circled branches in Fig.~\ref{fig_transition}F, which is nearly the iteration number on the two branches while much smaller for cells on the main trunk. Thus, the iterated mean first hitting time could distinguish cells at the branches or at the main trunk. If the ground truth is unknown, cells with large mean first hitting times as well as a significant gap between mean first hitting time to different fates (see the upper branch in Fig. \ref{fig_transition}G) could provide a reasonable candidate of the taboo set. This strategy works for the current synthetic example, however, its application to more practical examples needs more testing and deserves further study in the future.

\section{Discussion and conclusion}

The RNA velocity analysis has provided useful tools to predict future cell states within snapshot scRNA-seq data by modeling mRNA expression and splicing processes. Despite the established workflow and existing theoretical studies, several issues involved in parameter inference and downstream dynamical analysis of the current RNA velocity model remain elucidated, especially regarding the rationale and robust algorithm design and implementation. In this paper, we proposed several strategies to address these challenges through mathematical or statistical models as well as numerical analysis.

To unify the timescale of RNA velocity dynamics across different genes, we formulated the optimization framework based on either additive or multiplicative noise assumption to optimally determine the gene-specific rescaling parameters and proposed the numerical scheme to efficiently calculate the gene-shared latent time. Beyond the current independent gene assumption, it is possible to extend the current framework to RNA velocity models incorporating gene regulation or interactions\cite{bocci2022splicejac,wang2022inference,li2023multi,Zhang2023learning}.

To estimate the uncertainty of the inferred parameters and corresponding RNA velocity, we performed the confidence interval analysis of RNA velocity utilizing Fisher information approach, and SEM method\cite{MengUsing1991} was applied to obtain the asymptotic covariance of kinetic parameters in the dynamical RNA velocity model where latent cell time was involved in EM inference. In addition to the deterministic ODE model, our uncertainty quantification analysis and confidence interval construction approach could also be applied to stochastic RNA velocity model based on chemical reactions \cite{wang2022inference,Li1,gorin2022interpretable,gorin2022modeling}.

To determine the optimal hyper-parameters in the velocity-induced random walk, we analyzed its convergence rate toward continuous ODE dynamics in the operator sense, and assessed the dependence of convergence rate on sample size $n$ and data kernel bandwidth $\epsilon$. The results suggest that choosing kernel bandwidth $\epsilon$ around the scale of $O(n^{-2/(d+2)})$ provides the best operator approximation. It also indicates that as the dimensionality $d$ of the system increases, the accuracy of approximation would be impaired. Consequently, feature selection or dimensionality reduction could improve the cellular random walk construction. In parallel to the random walk approach, another strategy for downstream RNA velocity analysis is to fit the continuous velocity field \cite{qiu2022mapping} using vector-valued kernel methods \cite{ma2013regularized} as proposed in Dynamo \cite{qiu2022mapping} or neural-ODE methods\cite{chen2018neural,chen2022deepvelo,liu2022dynamical}, where similar convergence analysis could also be informative for the algorithm implementation.

To perform lineage inference and assign pseudo-time that is consistent with RNA velocity dynamics, we proposed to use the mean first hitting time of the velocity-induced random walk. The hitting time has been proposed for single-cell lineage analysis by defining lazy-teleporting random walk
on cellular similarity graph\cite{stassen2021generalized}. In the scVelo  package, a velocity pseudotime is defined based on the diffusion-like distance\cite{haghverdi2016diffusion,wolf2019paga} through the eigendecomposition of the weighted cellular velocity graph. In bifurcation systems with strong directionality induced by RNA velocity, our analysis and numerical examples suggest the introduction of taboo set for mean first hitting time analysis. Theoretically, the hitting time has close relation with the commute distance for undirected graph \cite{von2014hitting} $\mathcal{C}(x_i,x_j)=(T_i^j+T_ j^i)^\frac12$ where $T_i^j$ denotes the mean first hitting time from $x_i$ to $x_j$. It will be insightful to study the limit of mean first hitting time for the velocity-induced random walk when the sample tends to infinity, as studied for random geometric graph case \cite{von2014hitting} and unweighted directed graph case\cite{hashimoto2015random}.

Overall, the numerical analysis and statistical models presented in the current work could serve as a mathematical  step towards more robust and effective algorithmic implementation of the RNA velocity model  computation and analysis.

\appendix

\section{Proof of some lemmas} \label{appendixA}

\begin{lem}[Expansion of the un-normalized kernel $k_\epsilon$]\label{lem1}
The operator $\mathcal{G}_{\epsilon}$ for Gaussian-cosine scheme has the expansion
$$
\begin{aligned}
\mathcal{G}_{\epsilon} f(x) &=\frac{1}{\epsilon^{\frac{d}{2}}} \int k_{\epsilon}(x, y) f(y) \mathrm{d} y \\
&=m_{0} f(x)+\sqrt{\epsilon} m_{1} \mathcal{A} f(x)+O({\epsilon}),
\end{aligned}
$$
where
$$
\begin{aligned}
m_{0} &:=\mathcal{G}_{\epsilon} 1=\frac{1}{\epsilon^{\frac{d}{2}}} \int k_{\epsilon}(x, y) \mathrm{d}y \\
&=C_{d} \int_{0}^{\infty} r^{d-1} h\left(r^{2}\right) \mathrm{d} r \int_{-\pi}^{\pi}|\sin \theta|^{d-2} g(\cos \theta) \mathrm{d} \theta, \\
m_{1} &:=C_{d} \int_{0}^{\infty} r^{d} h\left(r^{2}\right) \mathrm{d} r \int_{-\pi}^{\pi} \cos \theta|\sin \theta|^{d-2} g(\cos \theta) \mathrm{d} \theta
\end{aligned}
$$
and
$$
\mathcal{A} f(x)=\|\nabla f(x)\| \cos \langle v(x), \nabla f(x)\rangle=\hat{v}(x) \cdot \nabla f(x), \hat{v}(x):=\frac{v}{\|v\|}.
$$
Here, $d>1, C_{d}=S_{d} / \int_{-\pi}^{\pi}|\sin \theta|^{d-2} \mathrm{~d} \theta$ and $S_{d}$ is the surface area of the $d$-dimensional unit sphere.
\end{lem}
The above lemma is in fact the Lemma 4.1 in \cite{Li1} except that the remainder term is explicitly characterized as $O(\epsilon)$ instead of $o(\sqrt{\epsilon})$. The proof is by straightforward derivations, which is also shown in Lemma \ref{lem2} below.

To get the variance error,  we first study the operator $\mathcal{\tilde{G}}_{\epsilon} $ defined by
$$
\mathcal{\tilde{G}}_{\epsilon} f^2(x)=\frac{1}{\epsilon^{\frac{d}{2}}} \int_{\mathbb{R}^{d}} (k_{\epsilon}(x, y) f(y) )^2q(y)\mathrm{d} y.
$$
The following lemma can be obtained.
\begin{lem}[Expansion of the kernel $k^2_\epsilon$]\label{lem2}
The operator $\mathcal{\tilde{G}}_{\epsilon}$ for Gaussian-cosine scheme has the expansion
$$
\begin{aligned}
\mathcal{\tilde{G}}_{\epsilon}f^2(x) &=\frac{1}{\epsilon^{\frac{d}{2}}} \int k_{\epsilon}^2 (x, y) f^2(y)q(y) \mathrm{d} y \\
&=\tilde{m}_{0} f^2(x)q(x)+\sqrt{\epsilon}\hat{m}_1\mathcal{A}(f^2(x)q(x))+O(\epsilon),
\end{aligned}
$$
where
$$
\tilde{m}_{0}
=C_{d} \int_{0}^{\infty} r^{d-1} h^2\left(r^{2}\right) \mathrm{d} r \int_{-\pi}^{\pi}|\sin \theta|^{d-2} g^2(\cos \theta) \mathrm{d} \theta
$$
and
$$
\tilde{m}_{1} =C_{d} \int_{0}^{\infty} r^{d} h^2\left(r^{2}\right) \mathrm{d} r \int_{-\pi}^{\pi} \cos \theta|\sin \theta|^{d-2} g^2(\cos \theta) \mathrm{d} \theta.
$$
Here, $d>1, C_{d}=S_{d} / \int_{-\pi}^{\pi}|\sin \theta|^{d-2} \mathrm{~d} \theta$ and $S_{d}$ is the surface area of the $d$-dimensional unit sphere.
\end{lem}

\begin{proof} For convenience,
let $v(x)=\|v\|(1,0, \cdots, 0)^{\mathrm{T}}$ without loss of generality.
Let us  first consider the case of $d=2$. Consider 2-dimensional  polar coordinates transformation
$$
\left\{\begin{array}{l}
y_{1}=x_{1}+r \cos \theta \\
y_{2}=x_{2}+r \sin \theta
\end{array}\right.
$$
where $\theta$ is the angle between $y-x$ and $v(x)$. Then we have
 \begin{equation}\label{v13}
\begin{aligned}
&\frac{1}{\epsilon}\int \left( k_{\epsilon}(x, y) f(y)\right)^2q(y) \mathrm{d} y\\
=&\frac{1}{\epsilon}\int_{0}^{\infty} \int_{-\pi}^{\pi} r
      h^2\left(\frac{r^{2}}{\epsilon}\right) g^2(\cos \theta) f^2(r, \theta)q (r, \theta)\mathrm{d} \theta \mathrm{d} r\\
 =&  \int_{0}^{\infty} r h^2\left(r^{2}\right) \int_{-\pi}^{\pi} g^2(\cos \theta)
        f^2(\sqrt{\epsilon} r, \theta)q (\sqrt{\epsilon}r, \theta) \mathrm{d} \theta \mathrm{d} r\\
 =&\int_{\epsilon^{\gamma-\frac{1}{2}}}^{\infty}+\int_{0}^{\epsilon^{\gamma-\frac{1}{2}}}\left(r h\left(r^{2}\right) \int_{-\pi}^{\pi} g^2(\cos \theta)
        f^2(\sqrt{\epsilon} r, \theta)q (\sqrt{\epsilon}r, \theta) \mathrm{d} \theta\right) \mathrm{d} r \\
:= &Q_{1}+Q_{2},
\end{aligned}
\end{equation}
where $0<\gamma<\frac{1}{2}$. Here
$$
\begin{aligned}
Q_{1} &=\int_{\epsilon^{\gamma-\frac{1}{2}}}^{\infty} r h\left(r^{2}\right) \int_{-\pi}^{\pi} g^2(\cos \theta) f^2(\sqrt{\epsilon} r, \theta)q (\sqrt{\epsilon}r, \theta)  \mathrm{d} \theta \mathrm{d} r \leq C \exp \left(-\epsilon^{2 \gamma-1}\right)=o(\epsilon).
\end{aligned}
$$
For $Q_{2}$, using Taylor expansion
$$
f^2(\sqrt{\epsilon} r, \theta)q(\sqrt{\epsilon} r, \theta)
=\left.f^2q\right|_{(0, \theta)}+\sqrt{\epsilon} r \Bigg(2fq\frac{\partial f}{\partial r}+ f^2\frac{\partial q}{\partial r}\Big|_{(0, \theta)}\Bigg)+O({\epsilon}),
$$
we get
\begin{align}
Q_{2}&=
    \int_{0}^{\epsilon^{\gamma-\frac{1}{2}}} r h^2\left(r^{2}\right) \int_{-\pi}^{\pi} g^2(\cos \theta) f^2(\sqrt{\epsilon} r, \theta)q(\sqrt{\epsilon} r, \theta) \mathrm{d} \theta \mathrm{d} r \nonumber\\
&= \int_{0}^{\epsilon^{\gamma-\frac{1}{2}}} r h^2\left(r^{2}\right) \int_{-\pi}^{\pi}
   g^2(\cos \theta)\left(\left.f^2q\right|_{(0, \theta)}+\left.\sqrt{\epsilon} r \left(2fq\frac{\partial f}{\partial r}+ f^2\frac{\partial q}{\partial r}\right|_{(0, \theta)}\right)+O(\epsilon)\right) \mathrm{d} \theta \mathrm{d} r \nonumber\\
&=\int_{0}^{\infty} r h\left(r^{2}\right) \int_{-\pi}^{\pi} g^2(\cos \theta)
    \left(\left.f^2q\right|_{(0, \theta)}+\left.\sqrt{\epsilon} r \left(2fq\frac{\partial f}{\partial r}+ f^2\frac{\partial q}{\partial r}\right|_{(0, \theta)}\right)\right) \mathrm{d} \theta \mathrm{d} r+O(\epsilon) \nonumber\\
&=\tilde{m}_{0} f^2(x)q(x)+\sqrt{\epsilon}\tilde{m}_1\mathcal{A}(f^2(x)q(x))+O(\epsilon). \label{v14}
\end{align}
For the high-dimensional case, the derivation is similar, so we omit it.
\end{proof}

\section{Proof of Theorem \ref{thm:discreteerror}} \label{appendixB}

\begin{proof}[Proof of Theorem \ref{thm:discreteerror}]
Following \cite{Singer}, we consider using the Chernoff inequality  to get an upper bound for $p(n,\alpha)$  with an $\alpha$-error. We will only estimate the term $P(\mathcal{L}_{\epsilon,n}f-\mathcal{L}_{\epsilon}f>\alpha)$ since the other part can be made similarly.

Let $\tilde{\alpha}=\sqrt{\epsilon}\alpha$. We have
$$
\begin{aligned}
p(n,\alpha) & =P\Big(\sqrt{\epsilon}(\mathcal{L}_{\epsilon,n}f-\mathcal{L}_{\epsilon}f)>\tilde{\alpha}\Big)\\
&=P\Bigg(\frac{\sum_{j=1}^{n}k_\epsilon(x,s_j)f(s_j)}{\sum_{j=1}^{n}k_\epsilon(x,s_j)}
-\frac{\int k_{\epsilon}(x, y)f(y) q(y) \mathrm{d} y}{\int k_{\epsilon}(x, y) q(y) \mathrm{d} y}>\tilde{\alpha}\Bigg).
\end{aligned}$$
Since $k_{\epsilon}(x,s_j)$ is positive, we have
$$
\begin{aligned}
p(n,\alpha)=P\Bigg(\sum_{j=1}^{n}\Big[\mathbb{E}\big(k_{\epsilon}(x,y)\big)k_{\epsilon}(x,s_j) f(s_j)
-\Big(\mathbb{E}\big(k_{\epsilon}(x,y)f(y)\big)+\tilde{\alpha} \mathbb{E}\big(k_{\epsilon}(x,y)\big)\Big)k_{\epsilon}(x,s_j)\Big]>0\Bigg),
\end{aligned}
$$
which is equivalent to
$$
\begin{aligned}
p(n,{\alpha})=P\Bigg(\sum_{j=1}^{n}Y_j>n\tilde{\alpha}\Big(\mathbb{E}\big(k_{\epsilon}(x,y)\big)\Big)^2\Bigg),
\end{aligned}
$$
 where
 \begin{equation}
\begin{aligned}
Y_j:=\Big[\mathbb{E}(k_{\epsilon}(x,y)) k_{\epsilon}(x,s_j) f(s_j)
-\mathbb{E}( & k_{\epsilon}(x,y) f(y)) k_{\epsilon}(x,s_j)\Big]\ +\\
&\tilde{\alpha} \mathbb{E}(k_{\epsilon}(x,y))\Big(\mathbb{E}(k_{\epsilon}(x,y))-k_{\epsilon}(x,s_j)\Big).
\end{aligned}
\end{equation}
We remark that the expectation $\mathbb{E}$ in the above and continued expressions are taken with respect to  the variable $y$ or $s_j$ whose probability density function is $q (y)$.

 It is easy to find that $Y_j$ are i.i.d random variables with $\mathbb{E}(Y_j)=0$.
 Next we calculate the variance of $Y_j$,
  \begin{equation}\label{v1}
\begin{aligned}
\mathbb{E}Y_j^2=K_1+K_2+K_3,
\end{aligned}
\end{equation}
where
$$
\begin{aligned}
K_1=&\left(\mathbb{E}(k_{\epsilon}(x,y))\right)^2\mathbb{E}(k_{\epsilon}^2(x,y)f^2(y)) -2\mathbb{E}(k_{\epsilon}(x,y))\mathbb{E}(k_{\epsilon}(x,y)f(y))\mathbb{E}(k^2_{\epsilon}(x,y)f(y))\\
&+\left(\mathbb{E}(k_{\epsilon}(x,y)f(y))\right)^2\mathbb{E}(k^2_{\epsilon}(x,y)),\\
K_2=&2\tilde{\alpha}\mathbb{E}(k_{\epsilon}(x,y))\left[\mathbb{E}(k_{\epsilon}(x,y)f(y))\mathbb{E}(k^2_{\epsilon}(x,y))
-\mathbb{E}(k^2_{\epsilon}(x,y)f(y))\mathbb{E}(k_{\epsilon}(x,y))\right],\\
K_3=&\tilde{\alpha}^2\left(\mathbb{E}(k_{\epsilon}(x,y))\right)^2\left[\mathbb{E}(k^2_{\epsilon}(x,y))-\left(\mathbb{E}(k_{\epsilon}(x,y))\right)^2\right].
\end{aligned}
$$
From Lemma \ref{lem1}, the expectations of $k_{\epsilon}(x,y)$ and $k_{\epsilon}(x,y)f(y)$ can be obtained
 \begin{equation}\label{v2}
\begin{aligned}
\mathbb{E}(k_{\epsilon}(x,y))&=\epsilon^{\frac{d}{2}}\left(m_0q(x)+\sqrt{\epsilon}m_1\hat{v}(x)\cdot \nabla q(x)+O(\epsilon)\right),\\
\mathbb{E}(k_{\epsilon}(x,y)f(y))&=\epsilon^{\frac{d}{2}}\left(m_0q(x)f(x)+\sqrt{\epsilon}m_1\hat{v}(x)\cdot \nabla (q(x)f(x))+O(\epsilon)\right).
\end{aligned}
\end{equation}
From Lemma \ref{lem2}, the second moments $\mathbb{E}(k^2_{\epsilon}(x,y))$, $\mathbb{E}(k^2_{\epsilon}(x,y)f(y))$ and $\mathbb{E}(k^2_{\epsilon}(x,y)f^2(y))$ can be obtained
 \begin{equation}\label{v3}
\begin{aligned}
\mathbb{E}(k^2_{\epsilon}(x,y))&=\epsilon^{\frac{d}{2}}\left(\tilde{m}_0q(x)+\sqrt{\epsilon}\tilde{m}_1\hat{v}(x)\cdot \nabla q(x)+O(\epsilon)\right),\\
\mathbb{E}(k^2_{\epsilon}(x,y)f(y))&=\epsilon^{\frac{d}{2}}\left(\tilde{m}_0q(x)f(x)+\sqrt{\epsilon}\tilde{m}_1\hat{v}(x)\cdot \nabla (q(x)f(x))+O(\epsilon)\right),\\
\mathbb{E}(k^2_{\epsilon}(x,y)f^2(y))&=\epsilon^{\frac{d}{2}}\left(\tilde{m}_0q(x)f^2(x)+\sqrt{\epsilon}\tilde{m}_1\hat{v}(x)\cdot \nabla (q(x)f^2(x))+O(\epsilon)\right).
\end{aligned}
\end{equation}
Substituting (\ref{v2}), (\ref{v3}) into (\ref{v1}), we get
 \begin{equation}\label{v4}
\begin{aligned}
\mathbb{E}Y_j^2=&2\tilde{\alpha}\mathbb{E}(k_{\epsilon}(x,y))\left[\mathbb{E}(k_{\epsilon}(x,y)f(y))\mathbb{E}(k^2_{\epsilon}(x,y))-\mathbb{E}(k^2_{\epsilon}(x,y)f(y))\mathbb{E}(k_{\epsilon}(x,y))\right]\\
  &+\tilde{\alpha}^2\left(\mathbb{E}(k_{\epsilon}(x,y))\right)^2\left[\mathbb{E}(k^2_{\epsilon}(x,y))-\left(\mathbb{E}(k_{\epsilon}(x,y))\right)^2\right]+C\epsilon^{\frac{3d+2}{2}},
\end{aligned}
\end{equation}
where $C$ is a constant.

Note that our interested regime is $\tilde{\alpha} \lesssim O(\epsilon)$ since we are estimating an $O(\sqrt{\epsilon})$ quantity, namely $\sqrt{\epsilon}\mathcal{L} f(x)$, so an error $\tilde{\alpha}$ larger than the estimated quantity is meaningless.  For the $\tilde{\alpha}$ term in (\ref{v4}), straightforward calculations based on \eqref{v2}-\eqref{v3} show
$$
\begin{aligned}
\mathbb{E}(k_{\epsilon}(x,y))\left[\mathbb{E}(k_{\epsilon}(x,y)f(y))\mathbb{E}(k^2_{\epsilon}(x,y))-\mathbb{E}(k^2_{\epsilon}(x,y)f(y))\mathbb{E}(k_{\epsilon}(x,y))\right]
=O(\epsilon^{\frac{3d+1}{2}}).
\end{aligned}
$$
Thus, the $\tilde{\alpha}$ and $\tilde{\alpha}^{2}$ terms of the variance (\ref{v4}) are negligible, and we obtain
$$
\mathbb{E} Y_{j}^{2}=C_1\epsilon^{\frac{3d+2}{2}},
$$
where $C_1$ is a constant. By the Chernoff inequality, we obtain
\begin{equation}\label{v5}
\begin{aligned}
p(n,{\alpha})\leq 2\mathrm{exp}\Bigg(-\frac{n\tilde{\alpha}^2\epsilon^{\frac{d}{2}}}{C_1\epsilon }\Bigg)=2\mathrm{exp}\Bigg(-\frac{n\alpha^2\epsilon^{\frac{d}{2}}}{C_1 }\Bigg).
\end{aligned}
\end{equation}
The inequality (\ref{v5}) means that the correct magnitude of $\alpha$ should be made such that
\begin{equation}\label{v6}
n\alpha^2\epsilon^{\frac{d}{2}}=O(1),
\end{equation}
i.e., $\alpha\sim O(n^{-\frac12}\epsilon^{-\frac{d}{4}})$.
\end{proof}




\normalem
\bibliographystyle{siam}
\bibliography{RNA}
\end{document}